\documentclass[12pt]{article}
\usepackage[utf8x]{inputenc}
\usepackage{amsfonts, amssymb, euscript,graphics}
\usepackage{amsmath}
\usepackage{amscd}
\usepackage{graphicx}
\usepackage{eepic}
\usepackage{comment}
\usepackage{psfrag}
\usepackage{wrapfig}

\newcommand{\R}{{\mathbb R}}
\newcommand{\E}{{\mathbb E}}

\newcommand{\Z}{{\mathbb Z}}

\newcommand{\<}{\langle}
\renewcommand{\>}{\rangle}
\renewcommand{\le}{\leqslant}
\renewcommand{\ge}{\geqslant}

\newcommand{\ind}{\mathop{\rm ind}}

\renewcommand{\d}{\partial}

\newtheorem{theorem}{Theorem}[section]
\newtheorem{lemma}[theorem]{Lemma}
\newtheorem{proposition}[theorem]{Proposition}
\newtheorem{statement}[theorem]{Statement}
\newtheorem{corollary}[theorem]{Corollary}
\newfont{\cmbb}{cmex10 at 7pt}

\newcommand*{\lo}[1]{{\raisebox{-.15ex}{\mathsurround=0pt\makebox{$\scriptstyle\mskip-.5mu#1$}}}}

\newcommand{\ign}[1]{{}}
\def \.{\mskip 1mu}
\def \?{\mskip -1mu}

\def\Top{\rm{Top}}
\def\Add{\rm{Add}}

\newcommand*{\AutT}{\mathrm{Aut}\mskip.5mu_{\mathrm T}}
\newcommand*{\Aut}{\mathrm{Aut}}

\renewcommand*{\d}{\partial}

\newcommand{\M}{{\mathcal M}}

\newcommand*\proof[1]{{\noindent\it Proof#1.\/\enspace\ignorespaces}}
\newcommand*\pr{\proof{}}

\def\sq{\raisebox{-.1ex}{\hbox{$\square$}}}

\def\m#1{\includegraphics{#1.eps}}

\def\qed{{\unskip\nobreak\hfill\hskip.5em\nobreak\hfil\sq
\parfillskip=0pt\medbreak}}

\makeatletter

\renewcommand{\abovecaptionskip}{.8ex}

\newcommand*{\capdelimiter}{.}
\renewcommand{\@makecaption}[2]{%
\vspace{\abovecaptionskip}%
\sbox{\@tempboxa}{#1\capdelimiter#2} \ifdim \wd\@tempboxa >\hsize
   #1\capdelimiter#2 \par
\else
 \global\@minipagefalse \hbox to \hsize {\hfil #1\capdelimiter#2 \hfil}%
\fi \vspace{\belowcaptionskip}}

\newcommand{\subs}[1]
   {\refstepcounter{subsection}
   \medskip\noindent
   {\it\arabic{section}.\arabic{subsection}.\hspace{.25em}\ignorespaces#1.}
          \ignorespaces}
\newcommand{\subsempty}
   {\refstepcounter{subsection}
   \medskip\noindent
   {\it\arabic{section}.\arabic{subsection}.\hspace{.25em}}}
   \@addtoreset{equation}{section}

\makeatother

\title{Morse theory on manifolds with boundary I. Strong Morse function, cellular structures and algebraic simplification of cellular differential.}

\author{Petr E. Pushkar
\thanks{National Research University Higher School of Economics, Russian Federation.
This work has been funded by the  Russian Foundation for
Basic Research under the Grants RFBR 18-01-00461a. Supported in part by the Simons Foundation}}

\date{}

\date{}

\begin{document}
\maketitle
\begin{abstract}
Main subject of the paper is a (strong) Morse function on a compact manifold with boundary. We construct a cellular structure and discuss its algebraic properties in this paper. Also we get an estimation on Arnold's question on a number of critical points of a Morse function with given boundary condition.
\end{abstract}

\section*{Introduction}


\subsempty \label{0} Classical Morse theory studies the relationship between the set of critical points of a
Morse function on a manifold and the topology of the manifold.
In this paper we consider the case of a manifold with boundary.
%

We consider so-called strong Morse functions.
Recall that a function $F$ defined on a compact manifold $M$ with boundary $\d M$ is called a Morse function
if
\begin{enumerate}
\item{all of its critical points are non-degenerate and are contained in the interior of~$M$;} \item{the
restriction $F|_{\lo{\d M}}$ is Morse function on the closed manifold $\d M$.}
\end{enumerate}
Denote by ${\rm Crit}(F)$ the set of all critical points of the function $F$. A Morse function $F$ is called
a {\it strong Morse function} if for any $x,y\in{\rm Crit}(F)\cup{\rm Crit}(F|_\lo{\d M})$ we have $F(x)\ne
F(y)$. We will refer to a germ of a strong Morse function along the boundary as a {\it strong Morse germ}.

Let $F$ be a strong Morse function on $M$, and let $\E$ be a field. The main result of this paper is the
existence of powerful combinatorial structure on ${\rm Crit}(F)$.
%
%
\ign{
Let $f$ be a strong Morse germ along the boundary $\d M$ of a manifold~$M$. We will construct a finite set
$\mathcal{P}_{\E}(f,M)$ of polynomials with nonnegative integer coefficients depending on $f$, $(M,\d M)$ and
$\E$ only. The crucial property of $\mathcal{P}_{\E}(f,M)$ is that for any strong Morse function $F$
extending the germ $f$, the polynomial $P_{\E}(F)(t)$ belongs to the set $\mathcal{P}_{\E}(f,M)$.
Thus we have the following theorem:

\begin{theorem}[Generalized Morse inequalities]
\label{polynomial-estimate} Suppose that $f$ is a strong Morse germ. Each strong Morse function $F\colon
M\to\R$ extending $f$ determines a polynomial $P_{\E}(F) \in \mathcal{P}_{\E}(f,M)$ and a polynomial
$K_{\E}(F)$ with nonnegative integer coefficients so that
$$
P(F)(t)=P_{\E}(F)(t)+(1+t)K_{\E}(F)(t).
$$
\end{theorem}

The set $\mathcal{P}_{\E}(f,M)$ is explicitly constructed (see Sec.~\ref{polyndef}, \ref{G-vertexes},
\ref{G-edges} and \ref{G-matching}) from topological data described in sec.~\ref{data} below.

Concerning the problem of estimating the number of critical points of an extension of a germ along the
boundary we obtain the following theorem, which is a corollary of Theorem~\ref{polynomial-estimate}.

\begin{theorem}
\label{weak Morse estimate} {\rm(\,Generalized weak Morse estimates)} Suppose $f$ is a strong Morse germ and
let $F$ be a Morse function extending~$f$. Then,

(1) the number of critical points of the function $F$ is greater than or equal to
$\min\limits_{P\in\mathcal{P}_{\E}(f,M)}P(1)$;

(2) the number of critical points of the function $F$ of index $i$ is greater than or equal to the minimal
coefficient of~$\,t^i$ among all the coefficients of~$\,t^i$ of polynomials in the set
$\mathcal{P}_{\E}(f,M)$.
\end{theorem}

\subs{Relation with the classical Morse inequalities}\label{Classical Morse}
%
We will briefly recall the celebrated Morse inequalities for manifolds with boundary \cite{Morse}. Consider a
germ $f$ of a function along the boundary of a manifold $M$. We say that a critical point of the function
$f|_\lo{\d M}$ is {\it outward directed} (respectively, {\it inward directed}) if the derivative of $f$ in
the direction of the outer normal to the manifold at this point is positive (respectively, negative). For a
function $F$ on $M$ the outward (respectively, inward) directed critical point of $F|_\lo{\d M}$  are defined
as those of the germ of $F$ along $\d M$. We denote the number of critical points of the function $F$ of
index $i$ by $m_i(F)$ and the number of inward directed critical points of the function $F|_\lo{\d M}$ by
$m_i^{\d}(F,M)$.

For any field~$\E$ and Morse function $F$ the numbers $M_i=m_i(F)+m_i^{\d}(F,M)$ and the numbers
$b_i^{\E}(M)=\dim H_i(M;\E)$ satisfy the Morse inequalities \cite{Morse} for a manifold with boundary
\begin{equation}\label{Morseinequality}
M_k-M_{k-1}+...+(-1)^kM_0\ge b_k^{\E}(M)-b_{k-1}^{\E}(M)+...+(-1)^{k}b_0^{\E}(M)
\end{equation}
where $k\in\{0,...,n-1\}$, as well as
$$
M_n-M_{n-1}+...+(-1)^nM_0= b_n^{\E}(M)-b_{n-1}^{\E}(M)+...+(-1)^{n}b_0^{\E}(M).
$$
It is known that this system of inequalities is equivalent to the following statement, similar to the
statement of Theorem~\ref{polynomial-estimate}. Namely, there exists a polynomial $k_{\E}(F)$ having
nonnegative coefficients such that
$$
\sum M_it^i=P_{\E}(M)(t)+(1+t)k_{\E}(F)(t),
$$
where $P_{\E}(M)(t)=\sum b_i^{\E}t^i$ is Poincar\'{e} polynomial of the manifold $M$. This equality is
equivalent to the following equality
$$
P(F)(t)=P_{\E}(M)(t)-P_{-}(F,M)(t)+(1+t)k_{\E}(F)(t),
$$
where $P_{-}(F,M)(t)=\sum m_i^{\d}(F,M)t^i$. Note that the polynomial $P_{\E}(M)(t)-P_{-}(F,M)(t)$ depends
only on the manifold $M$ and the germ $f$ of $F$ along $\d M$, we will denote $P_{-}(F,M)$ by $P_{-}(f,M)$.

We will show (see Sec.~\ref{compare}) that, given a strong Morse germ $f$ and a polynomial $P\in
\mathcal{P}_{\E}(f,M)$
there exists a polynomial $Q$ with nonnegative integer coefficients such that
$P(t)=P_{\E}(M)(t)-P_{-}(f,M)(t)+(1+t)Q(t)$. Thus the classical Morse inequalities are consequences of the
Theorem~\ref{polynomial-estimate}. For $M$ closed, we show in Sec.~\ref{ex1} that the set
$\mathcal{P}_{\E}(f,M)$ consists of a single element, namely, the Poincar\'{e} polynomial $P_{\E}(M)$ of~$M$.
Hence in the case of closed manifolds Theorem~\ref{polynomial-estimate} is equivalent to the classical Morse
inequalities.

\subs{Simplest example} \label{example} For certain important classes of germs the classical Morse
inequalities yield adequate estimates for Arnold's problem.
More often, however, the estimate derived from classical Morse equalities are weak or even vacuous.

\begin{wrapfigure}{r}{75pt}
\begin{center}
\psfrag{R}[][][1]{\footnotesize$\R$}
\psfrag{F}[][][1]{\footnotesize$F$}
\includegraphics[width=.2\textwidth]{interval.eps}
 \end{center}
 \caption{}
\label{interval}
\end{wrapfigure}
Consider, for example, the simplest manifold with boundary --- a closed interval --- and a function $F$ shown
on Fig.~\ref{interval}. Let $f$ be the germ of $F$ along the boundary. It is clear that any extension of $f$
will have at least two critical points. However, the polynomial $P(f)=P_{\E}(M)-P_{-}(F,M)$ is equal to zero.
Hence, the classical Morse inequalities estimate the number of critical points of a generic extension of $f$
from below by $0$. We show (see Sec.~\ref{ex1}) that the set $\mathcal{P}_{\E}(f,M)$ consists of a single
polynomial, $\mathcal{P}_{\E}(f,M)=\{1+t\}$, which guaranties at least two critical points by
Theorem~\ref{weak Morse estimate}.

\subs{Topological data} \label{data} The set $\mathcal{P}_{\E}(f,M)$  constructed in Sec.~\ref{polyndef} by
an explicit procedure  starting from the following data (all the homologies are counted with coefficients in
$\E$):

(1) The dimensions of the homologies $H_k(M)$, $H_k(M,\d M)$, and $H_k(\d M)$ for any $k$.

(2) The critical values, indices and types (inward or outward) of critical points of the function $f|_{\lo{\d
M}}$.

Let $c_1<...<c_N$ be all the critical values of the function $f|_{\lo{\d M}}$. We fix a choice of numbers
$a_1,...,a_{N+1}$, such that $a_1<c_1<a_2<...<a_{N+1}$.

(3) For any pair $i,j$, such that $1\le i < j \le N+1$, and for any $k$ the dimension of the $k$-th homology
of the pair $(\{f|_{\lo{\d M}}\le a_j\}, \{f|_{\lo{\d M}}\le a_i\})$.

(4) For any $j \in \{1,...,N+1\}$ and nonnegative $k$ the dimension of the subspace
$$\iota_*(H_k(\{f|_{\lo{\d M}}\le a_i\}))\cap \d^*H_{k+1}(M,\d M)\subset H_k(\d M)$$
is known. The mapping $\iota_*$ is induced by the natural inclusion $\{f|_{\lo{\d M}}\le a_i\}\hookrightarrow
\d M$ and $\d^*$ is the connecting homomorphism from an exact subsequence of the pair.

The construction of $\mathcal{P}_{\E}(f,M)$ is independent of the other parts of the paper.

\subs{Plan of the paper} In Section $1$ we consider an Arnold's example, generalizing example of
Sec.~\ref{example}, of a germ of a function along the boundary of $n$-dimensional ball and explain methods
and techniques of the paper. In Section $2$ we correspond to a strong Morse function an algebraic object ---
a pair of chain complexes with a preferred basis. This pair is defined up to some equivalence.???

\subs{Previous results} The problem of finding the condition under which a germ of a function along the
boundary can be extended into the interior without critical points was considered in ~\cite{BlankLaudenbach}
and \cite{Carlos}.

The Arnold's problem for a closed $n$-dimensional ball was considered by Barannikov~\cite{Serega}. In this
case the results of \cite{Serega} imply estimates which coincide with ours.

BLAGODARNOSTI?????
}

\ign{
\section{On results and techniques of the paper}

\subs{Arnold's problem and classical Morse inequalities} One can show that the algebraic number of critical
points of a generic extension of a strong Morse germ $f$ to $M$ is independent of the extension and equals to
$\chi(M)-\sum_i (-1)^im_i(f)$, where $m_i(f)$ is the number of inward directed critical points of $f|_\lo{\d
M}$ of index $i$.
Denote the number $\chi(M)-\sum_i (-1)^im_i(f)$ by $\chi(M,f)$. The absolute value of $\chi(M,f)$ gives a
rough estimate for Arnold's problem. This estimate could be deduced from the classical Morse inequalities
which give, in general, a stronger estimate in Arnold's problem. We also note here, that Arnold's question of
estimation of a number of critical points from below is vacuous without boundary conditions. Indeed, for a
connected compact manifold $M$ of dimension $n>1$, such that $\d M \ne \varnothing$, and any collection
$m_0,...,m_n$ of non-negative integers one can construct a strong Morse function $F$ on $M$ such that
$m_i(F)=m_i$.

The example considered in \ref{example} is a particular case of a more general construction  due to Arnold.
For a given constant $C$ and $n>1$ the construction produces a germ $h$ along the sphere $S^n$ bounding a
closed ball $B=B^{n+1}$ such that $\chi(B,h)=0$, and the number of critical points of a generic extension of
$h$ to the ball is at least $C$.
The construction starts from an auxiliary closed connected manifold $N$ of dimension $n+1$ such that $\sum
b_i^{\E}(N)\ge C$ and $\chi(N)=0$ and a function $H$ on $N$ having a finite number of critical points.
Consider an embedding $e\colon B \hookrightarrow N$ such that all critical points of $H$ are contained in the
interior of the image of the embedding. Denote by $h$ the germ of $e^*H$ along the sphere $S^n$. After a
slight perturbation of $e$ we can assume that $h$ is a strong Morse germ.

The germ $h$ has the properties desired . Indeed, let $F$ be a generic extension of $h$ to $B$.
Let $G$ denote the function on $N$ uniquely determined by $G|_{N\setminus e(B)}=H$ and $G\circ e = F$.
The function $G$ is a Morse function, hence, by Morse theory the number of critical points of $G$ is at least
$\sum  b_i^{\E}(N)$ ($\ge C$). By construction, the critical points of $G$ are contained in the interior of
$e(B)$. Therefore, the number of critical points of $F$ is equal to the number of critical points of $G$.
The number $\chi(B, h)$ is equal to the algebraic number of critical points of $G$, which equals to $\chi(N)$
($=0$) by Morse theory.

We show (see Sec.~\ref{vacuous}) that, if $N$ is a product of a closed manifold with a circle, the classical
Morse inequalities do not guarantee existence of critical points of a generic extension of $h$ to the ball.
At the same time our inequalities estimate the number of critical points from below by the sum of Betti
numbers of the manifold~$N$.

}

\subs{Pairs of complexes} The main algebraic object of the paper is a pair of chain complexes of vector spaces over $\E$ equipped with additional structure.
A pair of chain complexes arises in the following way.
Using Morse theory one associates to each strong Morse function $F$ on $M$  a pair $(X, Y)$ of
$CW$-complexes which is homotopy equivalent to the pair $(M,\d M)$. We note here, that the pair $(X, Y)$ is
not uniquely defined in general, it depends on choices of cell approximation in the construction.
It turns out that cells contained in $Y$ are in one-to-one correspondence with the critical points of
$F|_{\lo{\d M}}$.
Cells of $X$ which are not contained in $Y$ are in one-to-one correspondence with elements of the union of
the set of critical points of the function $F$ and the set of outward directed critical points of the
function $F|_\lo{\d M}$.
It turns out that there is a natural order on all cells of $X$, such that the cellular boundary of each cell
is either a linear combination of cells of lower order or zero.

Consider the pair of cellular chain complexes with coefficients in $\E$ of the pair $(X,Y)$. This is a pair
of graded vector spaces $(L_1,L_2)$ with the differential $\d$.
The pair $(L_1,L_2)$ is independent of $(X,Y)$ and has a preferred ordered basis which depends only on the
function $F$. The differential $\d$, in general, depends on the pair $(X, Y)$. The value of $\d$ on a basis
element is either a linear combination of basis elements having lower order or zero. We say that such a
differential is an $M$-differential.

The ambiguity in a choice of $(X,Y)$ leads to an arbitrariness of $M$-differentials. We refer to an upper
triangular group a group of all graded automorphisms of $L_1$ preserving $L_2$ and having upper triangular
matrices in the preferred basis. Upper triangular group acts on $M$-differentials by conjugation.

Let $\d$, $\d_1$ be two $M$-differentials on $(L_1,L_2)$ given by Morse theory. We show that there
exists an upper triangular automorphism $S$, such that $\d=S^{-1}\d_1S$.

\subs{Partition of $M$-differentials}By the consideration above, each strong Morse function $F$ corresponds
to an orbit $O_F$ of the action of the upper triangular group on the space of $M$-differentials.

We consider the space of all $M$-differentials acting on a pair of graded vector spaces equipped with an
ordered basis and the action by conjugation of the upper triangular group on this space.
Additional conditions motivated by topological reasons determine a subspace $\mathcal{D}$ of the space of all
$M$-differentials which is invariant under the action of the upper triangular group. 
Our main result is a partition of the set of orbits of this action on $\mathcal{D}$ into a finite number of subsets. 
This partition is generalization of Barannikov normal forms of Morse complexes of Morse functions on a closed manifold \cite{Serega}.

\subs{Morse type estimation} We apply that partition to get estimates in the following question of Arnold: what could be an estimation from below of a number of critical points of  Morse extension of a strong Morse germ along the boundary of a compact manifold?

We construct numbers $n_1(g,M,\E)$ and $n_2(g,M,\E)$ and estimate total number of a critical points of any Morse extension of a germ $g$ along $\d M$ by  $n_1(g,M,\E)$ and $n_2(g,M,\E)$ (see Theorems~\ref{n_1} and \ref{n_2}).  Also we give a criterium of existence of critical points (possibly degenerate) of any extension of $g$ to $M$ in a simple Theorem~\ref{criterium}.

\ign{
We
show that each subset in the resulting partition has a canonical representative which decomposes into a
direct sum of differentials of sixteen different types.

Thus, it turns out that there is a remarkable combinatorial structure on the set ${\rm Crit}(F)\cup{\rm
Crit}(F|_\lo{\d M})$: critical points of the function and of its restriction to the boundary can be naturally
divided into sets (consisting of one, two, three, and four elements) of the sixteen different types. This
combinatorial structure is a generalization to the case of a manifold with boundary of the division of the
critical points of a strong Morse function on a closed manifold into pairs and points ``responsible for
homologies''. In particular, this combinatorial structure gives rise to the partition ${\rm
Crit}(F)=\Top_\E(F)\cup\Add_\E(F)$ mentioned in sec.~\ref{0}.

In addition, we associate a finite graph $\Gamma(O)$ to an orbit $O$ of the action of the upper triangular
group on $\mathcal{D}$.
Let $f$ be a strong Morse germ along the boundary $\d M$ of a manifold $M$ and let $F$ is a strong Morse
function on $M$ extending the germ $f$. We show that the graph $\Gamma(O_F)$ depends only on the manifold $M$
and the germ $f$, $\Gamma(O_F)=\Gamma_\E(M,f)$. The topological data needed for construction of
$\Gamma_\E(M,f)$ is described below???.

??? Skazat' (snova - proverit' snowa li eto) chto mnogohleny stroqtsq po grafu

Denote the number of vertices of the graph $\Gamma$ by $v(\Gamma)$. Recall that a matching is a collection of
edges without common vertices. Let $m(\Gamma)$ equals two times the maximal number of edges in a matching.

We show that the number $\min\limits_{P\in\mathcal{P}_{\E}(M,f)}P(1)$ from Theorem~\ref{weak Morse estimate}
has the following interpretation in terms of $\Gamma_\E(M,f)$.

\begin{theorem}
\label{estimate} The number of critical points of a Morse function continuing a strong Morse germ $f$ is
greater then or equal to
$$
v(\Gamma_\E(M,f)) - m(\Gamma_\E(M,f))=\min\limits_{P\in\mathcal{P}_{\E}(M,f)}P(1).
$$
\end{theorem}

The set of vertexes of $\Gamma_\E(M,f)$ are, by definition, graded by integers. A $k$-th component of the set
of vertexes of $\Gamma_\E(M,f)$ is, by construction, the disjoint union of five sets $\mathrm{A}_k$,
$\mathrm{B}_k$, $\mathrm{C}_k$, $\mathrm{D}_k$ and $\mathrm{E}_k$.
The following theorem generalize (in terms of $\Gamma_\E(M,f)$) weak Morse inequalities $m_k(F)\ge
b_k^{\E}(M)-m_k(f,\d M)$ ($m_k(f,\d M)$ is the number of index $k$ inward critical points of $f$).

\begin{theorem}\label{weak Morse-graph}
The number $m_i(F)$ of critical points of index $k$ of a Morse function continuing a strong Morse germ $f$ is
greater then or equal to
$$
\#\mathrm{B}_k+\#\mathrm{C}_k+\#\mathrm{D}_k+\#\mathrm{E}_k-\#\mathrm{A}_{k-1}.
$$
\end{theorem}
}

\section{Functions on manifolds with boundary and pairs of complexes}

The standard procedure of Morse theory \cite{Milnor} associates a $CW$-complex to a Morse function on a
closed manifold.
Starting from a strong Morse function on a manifold $M$ with the boundary $\d M$ we construct a pair $(X,Y)$
of $CW$-complexes which is homotopy equivalent to the pair $(M,\d M)$.
In general, the pair $(X,Y)$ is not uniquely defined .
It depends on cellular approximations used in the construction below.
We study (at the level of cellular differentials) the ambiguity in our construction.

\subs{Bifurcations of sublevel sets} \label{CW}
Let $F$ be a strong Morse function on a~manifold $M$ with the boundary $\d M$.
We denote a sublevel set $\{F\le c\}$ by $F_c$, and the set $\{F|_{\lo{\d M}}\le c\}$ by $F^{\d}_c$.
Let $c_1<...<c_N$ be critical values of the functions $F$ and $F|_{\lo{\d M}}$.
For a topological space $X$ and a continuous map $\varphi\colon S^{k-1}\to X$ we denote by $X\cup_\varphi
e^k$ the result of attaching a cell $e^k$ of dimension $k$ along $\varphi$ to $X$.
Recall that a pair of topological spaces $(A,B)$ is a strong deformation retract of a pair $(A_1,B_1)$, if
$(A_1,B_1)\supset (A,B)$ and there exists a family $f_{t, t\in[0,1]}\colon A_1 \to A_1$ of continuous maps
such that $f_0=Id$, $f_t(B_1)\subset B_1$,  $f_t|_A=Id$ for any $t\in [0,1]$ and $f_1(A_1)=A$, $f_1(B_1)=B$.

The topology of the pair $(F_c,F^{\d}_c)$ changes when the parameter $c$ goes through critical values as
follows:
\begin{proposition}
\label{bifurcation} (0) If an interval $[a,b]$ does not contain critical values $c_1,...,c_N$, then the pair
$(F_a,F^{\d}_a)$ is a strong deformation retract of the pair $(F_b,F^{\d}_b)$.

Take  $c\in\{c_1,...,c_N\}$ and a sufficiently small number $\varepsilon>0$.
Consider the pairs $(F_{c-{\varepsilon}},F^{\d}_{c-{\varepsilon}}) \subset (F_{c-{\varepsilon}}\cup
F^{\d}_{c+\varepsilon},F^{\d}_{c+\varepsilon})\subset (F_{c+\varepsilon},F^{\d}_{c+\varepsilon})$.

(1) Let $c$ be the value of the function $F$ at a critical point of index $k$.
The pair $(F_{c-{\varepsilon}},F^{\d}_{c-{\varepsilon}})$ is a strong deformation retract of
$(F_{c-{\varepsilon}}\cup F^{\d}_{c+\varepsilon},F^{\d}_{c+\varepsilon})$. There exist an attaching map
$\varphi$ and a homotopy equivalence
$$
h\colon (F_{c+\varepsilon},F^{\d}_{c+\varepsilon})\to ((F_{c-{\varepsilon}}\cup
F^{\d}_{c+\varepsilon})\cup_\varphi e^k,F^{\d}_{c+\varepsilon}),
$$
which is the identity on $F_{c-{\varepsilon}}\cup F^{\d}_{c+\varepsilon}$.

(2) Let $c$ be the value of the function $F|_{\lo{\d M}}$ at an inward directed critical point of index $k$.
There exist an attaching map $\varphi$ and a homotopy equivalence
$$
h\colon (F_{c-{\varepsilon}}\cup F^{\d}_{c+\varepsilon},F^{\d}_{c+\varepsilon})\to
(F_{c-{\varepsilon}}\cup_\varphi e^k,F^{\d}_{c-{\varepsilon}}\cup_\varphi e^k),
$$
which is the identity on $F_{c-{\varepsilon}}$.
The pair $(F_{c-{\varepsilon}}\cup F^{\d}_{c+\varepsilon},F^{\d}_{c+{\varepsilon}})$ is a strong deformation
retract of $(F_{c+\varepsilon},F^{\d}_{c+\varepsilon})$.

(3) Let $c$ be the value of the function $F|_{\lo{\d M}}$ at an outward directed critical point of index $k$.
There exist an attaching map $\varphi$ and a homotopy equivalence
$$
h\colon (F_{c-{\varepsilon}}\cup F^{\d}_{c+\varepsilon},F^{\d}_{c+\varepsilon}) \to
(F_{c-{\varepsilon}}\cup_\varphi e^k,F^{\d}_{c-{\varepsilon}}\cup_\varphi e^k),
$$
which is the identity on $F_{c-{\varepsilon}}$. Moreover, there exist an attaching map $\varphi_1$ and a
homotopy equivalence
$$
h_1\colon (F_{c+{\varepsilon}},F^{\d}_{c+{\varepsilon}})\to ((F_{c-{\varepsilon}}\cup
F^{\d}_{c+\varepsilon})\cup_{\varphi_1}e^{k+1},F^{\d}_{c+\varepsilon}),
$$
which is the identity on $F_{c-{\varepsilon}}\cup F^{\d}_{c+\varepsilon}$. The space $F_{c-{\varepsilon}}$ is
a strong deformation retract of the space $F_{c+\varepsilon}$. \qed
\end{proposition}
Proposition \ref{bifurcation} is a relative version of standard \cite{Milnor} results from Morse theory, see \cite{JR}.
Its proof is parallel to the standard considerations, and follows from the relative version of the Morse
lemma, which states that for each inward (respectively, outward) critical point of $F|_{\lo{\d M}}$ with the
critical value $c$ there exist local coordinates $(x,y)$ $(y \in \R_+)$ centered at the critical point such
that $F(x,y)=c + y +Q(x)$ (respectively, $F(x,y)=c - y +Q(x)$) where $Q$ is a ``sum of squares'', and from an
explicit description of cells and retractions for such coordinate choice. We omit the details.

\begin{figure}[htb]
\begin{center}
\includegraphics[width=5in]{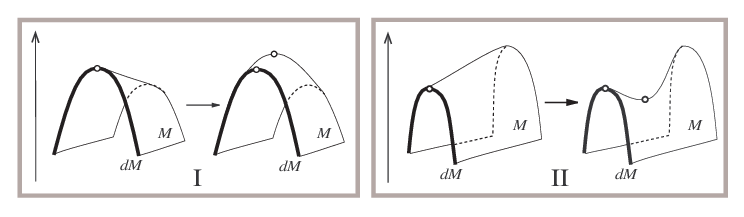}
\caption{Surgery on a strong Morse function}\label{Fig:Morsedeform}
\end{center}
\end{figure}

\subs{Remark} There exists a surgery on a strong Morse function which eliminates its outward (respectively,
inward) critical points and does not change the restriction to the boundary.
This surgery is local, that is defined in a collection of neighborhoods of critical points of the restriction
to the boundary. 
Two-dimensional examples of such a surgery are shown in Figure \ref{Fig:Morsedeform}. 
%
%
Each surgery adds an additional critical point inside the manifold.
The surgery eliminating all inward points was used in~\cite{Morse}. 
It is easy to obtain the classical Morse inequalities by combining this surgery with the results of parts (0),(1) and (3) of
Proposition~\ref{bifurcation}.
%

\subs{Morse chain}\label{constrCW}
We continue with the notations introduced in Section~\ref{CW}.
The function $F$ takes its maximal value at a point from the set ${\rm Crit}(F)\cup {\rm Crit}(F|_\lo{\d
M})$. Hence, the value set of the function $F$ belongs to the interval $[c_1,c_N]$.
We fix numbers $a_0,...,a_{N}$ such that $a_0<c_1<a_1<...<a_{N-1}<c_N<a_{N}$.
Consider the chain of inclusions of topological spaces:
\begin{multline*}
(\varnothing,\varnothing)=(F_{a_0},F^{\d}_{a_0})\subset (F_{a_0}\cup F^{\d}_{a_1}, F^{\d}_{a_1})\subset
(F_{a_1}, F^{\d}_{a_1})\subset...\\
...\subset(F_{a_{N-1}}, F^{\d}_{a_{N-1}})\subset(F_{a_{N-1}}\cup F^{\d}_{a_N},
F^{\d}_{a_N})\subset(F_{a_{N}}, F^{\d}_{a_{N}})=(M,\d M).
\end{multline*}
There are $2N+1$ pairs in the chain, we denote them as
$(U_0,V_0)\subset(U_1,V_1)\subset...\subset(U_{2N},V_{2N}).$

According to Proposition~\ref{bifurcation}, either $U_{i+1}$ is homotopy equivalent to $U_i$ or $U_{i+1}$ is
homotopy equivalent to $U_i$ with a cell attached.
Applying the standard technique (see \cite{Milnor}) one can construct a chain of inclusions of pairs of
$CW$-complexes
$$
(\widetilde{X}_0,\widetilde{Y}_0)\subset(\widetilde{X}_1,\widetilde{Y}_1)\subset...\subset(\widetilde{X}_{2N},\widetilde{Y}_{2N})$$
and homotopy equivalences $\widetilde{h}_i\colon(U_i,V_i)\to (\widetilde{X}_i,\widetilde{Y}_i)$ for
$i\in\{0,...,2N\}$ such that the diagram
$$
\begin{CD}
(U_0,V_0) & \.\.\.\subset \.\.\.& (U_1,V_1) & \.\.\.\subset\.\.\.& ... & \.\.\.\subset\.\.\.&(U_{2N},V_{2N})\\
@V{\widetilde{h}_0}VV  @V{\widetilde{h}_1}VV  ... && @V{\widetilde{h}_\lo{2N}}VV\\
(\widetilde{X}_0,\widetilde{Y}_0)&\.\.\.\subset\.\.\.&(\widetilde{X}_1,\widetilde{Y}_1)&\.\.\.\subset\.\.\.&
... &\.\.\.\subset\.\.\.&(\widetilde{X}_{2N},\widetilde{Y}_{2N})
\end{CD}
$$
is commutative and satisfies the following condition: for $i\in\{0,...,2N-1\}$,
$\widetilde{Y}_{i}=\widetilde{Y}_{2N}\cap \widetilde{X}_i$ and $\widetilde{X}_{i+1}$ is either equal to
$\widetilde{X}_{i}$ or is the result of attaching of a single cell to $\widetilde{X}_{i}$.

We recall standard topological notions. Let $(A_0,B_0)\subset(A_1,B_1)\subset...\subset(A_K,B_K)=(A,B)$ and
$(C_0,D_0)\subset(C_1,D_1)\subset...\subset(C_K,D_K)=(C,D)$ be filtered pairs of topological spaces. A
filtered (continuous) map is a map of pairs $h\colon(A,B)\to (C,D)$ such that $h(A_i)\subset C_i$,
$h(B_i)\subset D_i$ for any $i \in \{0,...,K\}$. A filtered homotopy between filtered maps $h_j:(A,B)\to
(C,D)$, $j\in\{0,1\}$ is a filtered map $H\colon(A\times I, B\times I)\to(C,D)$ such that
$H|_{{A\times\{j\}}}=h_j$, $j\in\{0,1\}$.
Two filtered maps $h_j:(A,B)\to (C,D)$, $j\in\{0,1\}$ are called {\it filtered homotopic} if there exists a
filtered homotopy between them.
A filtered map $h\colon(A,B)\to (C,D)$ is {\it a filtered homotopy equivalence} if there exists a filtered
map $g\colon(C,D)\to (A,B)$, such that the maps $h\circ g$ and $g\circ h$ are filtered homotopic to $Id_A,
Id_C$ respectively. It is easy to show along the lines of \cite{Milnor}, pp. 20--23, that the map
$\widetilde{h}_{2N}$ above is a filtered homotopy equivalence.

We say that {\it a Morse chain $\mathbf{M}$ of a strong Morse function~$F$} is a following triple:
\begin{enumerate}
\item{a $CW$-pair $(X,Y)$;}
\item a $CW$-filtration $(\varnothing,\varnothing)=(X_0,Y_0)\subset...\subset(X_{2N},Y_{2N})=(X,Y)$, such
that, for each $i\in\{0,...,2N-1\}$, ${Y}_{i}={Y}_{2N}\cap {X}_i$ and $X_i$ is either equal to $X_{i-1}$ or
is the result of attaching a single cell to $X_{i}$;
\item{ a filtered homotopy equivalence $h\colon (M,\d M)\to (X,Y)$.}
\end{enumerate}
We will assume below that orientation of cells in a Morse chain is somehow fixed.

In general, the complexes $X_i,Y_i$ from a Morse chain are not uniquely defined.
%
However, for any $i\in\{0,...,2N\}$ the number of cells of a given dimension in $X_{i},Y_{i}$ is determined
by the function~$F$ only.
According to Proposition~\ref{bifurcation}, the total number $T=T(F)$ of cells in the complex $X_{2N}$ is the
number of outward critical points of the function $F|_\lo{\d M}$ plus the number of critical points of the
function $F$. 
The total number of cells in the complex $Y_{2N}$ is equal to the number of critical points of
the function $F|_\lo{\d M}$.


\subs{Upper-triangularity} Consider a Morse chain $\mathbf{M}$ of a strong Morse function $F$.
We enumerate the cells of $\mathbf{M}$ by $e_1(\mathbf{M}),...,e_T(\mathbf{M})$ in the order of their
appearance in the subcomplexes $X_i$ --- a cell attached later has a bigger number than a cell attached on
earlier step.

Take two Morse chains $\mathbf{M}$, $\mathbf{M}'$ of the same strong Morse function~$F$. Let $\mathbf{M}$
consist of pairs $\{(X_i,Y_i)\}$ and a filtered homotopy equivalence $h$ and we mark the same objects for
$\mathbf{M}'$ with primes. Consider a filtered homotopy inverse
$g\colon (X_{2N},Y_{2N})\to(U_{2N},V_{2N})$ of the filtered map $h$.
Set $S=h\circ g$. The map $S$ induces a cellular homotopy equivalence $S_i\colon X_i\to X'_i$ for any
$i\in\{1,...,2N\}$. Let $S^\#_i$ denote the induced map of the cellular chain complexes.
The map $S^\#_{2N}$ is an isomorphism of complexes, moreover, the following statement holds.
\begin{proposition}\label{upper} The matrix of $S^\#_{2N}$ with respect to the bases
$(e_1(\mathbf{M}),...,e_L(\mathbf{M}))$ and $(e_1(\mathbf{M}'),...,e_L(\mathbf{M}'))$ is upper-triangular
with $\pm 1$ on the diagonal.
\end{proposition}
\pr Consider a cell $e_k(\mathbf{M})$. Consider the minimal $i=i(k)$ such that $e_k(\mathbf{M})\in X_{i}$. By
definition $X_i$ consists of $k$ cells $e_1(\mathbf{M}),...,e_k(\mathbf{M})$.
The map $S$ is filtered, hence $S^\#_{2N}(e_k(\mathbf{M}))=S^\#_{i}(e_k(\mathbf{M}))$. The number of cells in
$X'_i$ is equal to the number of cells in $X_i$. Therefore, $S^\#_{i}(e_k(\mathbf{M}))$ is a linear
combination of cells $e_1(\mathbf{M}'),...,e_k(\mathbf{M}')$.
This proves that the matrix of $S^\#_{2N}$ is upper-triangular.

A diagonal element $(S^\#_{2N})_{k,k}$ is the degree of the map
$$
S^{\dim e_k(\mathbf{M})}=X_{i(k)}/X_{i(k)-1}\to X'_{i(k)}/X'_{i(k)-1}=S^{\dim e_k(\mathbf{M}')}
$$ induced from $S_{i(k)}$.
Since $S_{i(k)}$ is a homotopy equivalence, this degree is equal to $\pm 1$.\qed

The considerations above motivate the following definitions.

\subs{Definition of $M$-complexes and of their isomorphisms}\label{Mdef} An {\it$M$-complex} is a following
structure:
\begin{enumerate}
\item A finite complex of finite-dimensional vector spaces over a field $\E$:
$$
0 \longrightarrow C_K \overset{\d_K}{\longrightarrow} C_{K-1} \overset{\d_{K-1}}{\longrightarrow} ...
\overset{\d_{L+1}}{\longrightarrow} C_L \overset{\d_L}{\longrightarrow} C_{L-1} \longrightarrow 0
$$
(that is, $\d_i\circ\d_{i+1}=0$ for all $i$).
We denote the direct sum $\oplus_i \d_i$ by~$\d$.
\item Every space $C_i$ is equipped with a fixed basis.
\item A linear order is chosen on the union $A$ of all bases, so that the ``decreasing order'' condition is
satisfied: for any $a\in A$ either the vector $\d(a)$ is a linear combination of elements from $A$ of orders
smaller then $a$ or zero.
\end{enumerate}

We will also need another equivalent definition of an $M$-complex. Let $A$ be a finite linearly ordered
graded set $\{a_1,...,a_N\}$, $a_1\prec...\prec a_N$.
A grading on $A$ is a mapping to $A\to\Z$, we denote it by ${\rm deg}$. The number ${\rm deg}(a)$ is called
the {\it degree} of the element $a\in A$.
The vector space $\E(A)=\E\otimes A$ of all formal linear combination of elements of $A$ with coefficients in $\E$  is naturally graded.
An {\it $M$-differential} ($M_A$-differential if we need to be more specific) is a differential $\d$ on
$E(A)$ of degree $-1$, such that $\d(\E\otimes\{a_1,...,a_i\})\subset \E\otimes\{a_1,...,a_{i-1}\}$
for all $i\in\{1,...,N\}$.
An $M$-complex is a graded vector space $\E(A)$ with a $M$-differential $\d$. We denote it by
$\M_{A,\? \d}$.
$M$-complexes equipped with additional structure were considered in \cite{Serega} under the name of framed
Morse complexes.

We say that two $M$-complexes $\M_{A_1, \d_1}$ and $\M_{A_2, \d_2}$ are equal, if the sets $A_1$ and $A_2$
are graded ordered isomorphic and the matrices of the differentials $\d_1$ and $\d_2$ in the bases $A_1$ and
$A_2$ coincide.
By $\Aut(A)$ we denote a group of all graded automorphisms of the space $\E(A)=\E\otimes A$. 
By $\AutT(A)\subset\Aut(A)$ we denote a subgroup of all graded automorphisms of the space $\E(A)$,
preserving each vector subspace $\E\otimes\{a_1,...,a_i\}$, $i\in\{1,...,N\}$.
The matrices in the basis~$A$ of operators from the group $\AutT(A)$ are upper-triangular.
We say that  $M_A$-differentials $\partial_1\?,\?\partial_2\colon \E(A) \to \E(A)$ are
equivalent  (or $A$-equivalent), if there exists $g \in \AutT(A)$, such that $\partial_{2} = g
\.\.\partial_{1}g^{-1}$. 
%
We say that two $M$-complexes are {\it isomorphic}, if they become equal after
replacing of one of the $M$-differentials by an equivalent.

\subs{Pairs of $M$-complexes} Let  $\M_{A,\? \d}$ be an $M$-complex.
For a subset $B$ of $A$ such that $\d(\E(B))\subset \E(B)$,
$\M_{B,\,\d|_\lo{\E(B)}}$ is an $M$-complex (the order and the grading on $B$ are induced from those
on $A$). We will say that $\M_{B,\,\d|_\lo{\E(B)}}$ is an $M$-subcomplex of the $M$-complex
$\M_{A,\d}$ and will denote $\d|_\lo{\E(B)}$ by $\d_{B}$.
A pair consisting of  an $M$-complex and its $M$-subcomplex will be called  a pair of $M$-complexes or an
$M$-pair.
The differential $\d$ in that case will be called an $M_{A,B}$-differential. We denote an $M$-pair $(\M_{A,\?
\d},\M_{B,\,\d_\lo{B}})$ by $\M_{A,B,\,\d}$.

By $\Aut(A,B)\subset \Aut(A)$ we denote a subgroup of all graded automorphisms $g\in \Aut(A)$ such
that $g(\E\otimes B)=\E\otimes B$.
Consider the subgroup $\AutT(A,B)$ of the group $\AutT(A)$, consisting of all elements $g\in\AutT(A)$ such
that $g(\E\otimes B)=\E\otimes B$. Clearly $\AutT(A,B)=\Aut(A,B)\cap\AutT(A)$. 
Two $M_{A,B}$-differentials $\d_1$ and $\d_2$ are called
$(A,B)$-equivalent (or equivalent), if $\d_2=g\d_1g^{-1}$ for some $g\in\AutT(A,B)$.
Two $M$-pairs $\M_{A_1,B_1,\,\d_1}$ and $\M_{A_2,B_2,\,\d_2}$ are called {\it isomorphic}, if they become
equal after changing of one of the $M$-differentials on an $(A_1,B_1)$-equivalent.

\subs{Algebraic model of a strong Morse function}\label{amodel} 
%
Let $F$ be a strong Morse function and $\mathbf{M}$ be a Morse chain of $F$.
The cellular boundary of $e_k(\mathbf{M})$ is either zero or linear combination of cells with smaller
indices.
Hence, a Morse chain naturally generates a pair of $M$-complexes $\M_{A_F,B_F,\,\d}$ with the
$M$-differential $\d=\d(\mathbf{M})$.
The set $B_F$ may be identified with critical points of the function $F|_{\lo{\d M}}$ graded by the Morse
index $\ind_M$ and ordered with respect to critical values.
The set $A_F$ is a result of the following operations. Firstly we add to $B_F$ the set of critical points of
function $F$ graded by~$\ind_M$. Denote the resulting set by $X_F$. 
It is naturally ordered with respect to the critical values. 
Denote by $C_F\subset B_F$ a subset consists of all outward critical points of the
function $F|_{\lo{\d M}}$. 
For each element $b\in C_F$ we add to $X_F$ an element $b_+$ next to $b$ with
degree $\ind_M F|_{\lo{\d M}}(b)+1$. 
The resulting set is $A_F$.

The following statement summarizes the previous observations.
\begin{statement}
\label{F-CW} A strong Morse function~$F$ naturally corresponds to a pair of $M$-complexes
$\M_{A_F,B_F,\,\d}$, defined up to an isomorphism. \qed
\end{statement}

An arbitrary pair of $M$-complexes $\M_{A_F,B_F,\,\d}$ isomorphic to an $M$-pair constructed following Morse
theory from a strong Morse function $F$ will be called an {\it algebraic model of function~$F$}.
The boundary of a $1$-cell consists of at most two $0$-cells. Hence, even if we consider integer coefficients
and isomorphisms only, then not every algebraic model corresponds to a Morse chain from
Section~\ref{constrCW}.

\begin{figure}[htb]
\begin{center}
\includegraphics[width=5 in]{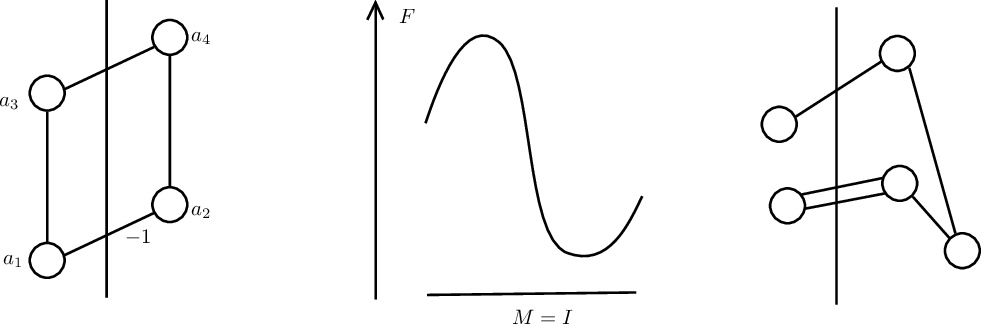}
\caption{\,\,Examples of $M$-pairs}\label{Fig:M-pairexmpls}
\end{center}
\end{figure}

\subs{Drawing of a pair of $M$-complexes} An $M$-pair $\M_{A,B,\d}$ we will draw as follows. 
Elements of the
set $A$ we draw by circles and place this circles along the vertical axis in correspondence with the order on
$A$: a circle corresponding to an element $a_i$ is higher then a circle corresponding to an element $a_j$ if
$i>j$. 
Circles corresponding to elements from the set $B$ we draw left to the vertical axis, circles
corresponding to elements from the set $A\setminus B$ we draw right to the vertical axis. 
If $\d a_i =\sum\limits_{k\in I} \lambda_k a_k$, $\lambda_k\ne 0$ then we connect circles corresponding to elements $a_i$and $a_k, k \in I$ by segments labelled by $\lambda_k$ if $\lambda_k\ne 1$. 
For example, an $M$-pair
$\M_{A,B,\d}$, where $A=\{a_1,a_2,a_3,a_4\}$, $B=\{a_1,a_3\}$, and differential $\d$ defined on basis
elements: $\d a_4=a_3+a_2$, $\d a_3 =a_1$, $\d a_2 =-a_1$, $\d a_1 = 0$ is shown on the left side of
Fig.~\ref{Fig:M-pairexmpls}. 

%
We show below that each outward critical point $b\in C_F$ appears with
non-zero coefficient in $\d(b_+)$. 
We will draw double segment connecting $b$ with $b_+$ instead of ordinary
segment. 
On the Fig.~\ref{Fig:M-pairexmpls}. (center and right) we show a graph of a function~$F$ on a segment and its algebraic
model with $\Z_2$~-coefficients.
In that case $A_F$ consists of five, $B_F$ of two and $C_F$ of one element.
%

\section{Pair of M-complexes. Formulation of the Main Theorem.} 
\label{formulation}

An arbitrary $M$-pair considered up to the equivalence is a difficult object to work with. 
We will stratify
the set of all $M$-pairs into pieces so that each stratum contains a relatively simple $M$-pair encoding the
stratum.

\subs{$\d$-trivial elements} Consider an $M$-pair $\M_{A,B,\,\d}$, $A=\{a_1\prec ...\prec a_N\}$.
We call an element $a_k \in B$ {\it $\d$-trivial}, if $k<N$,  $a_{k+1}\in A\setminus B$ and $a_{k}$ appears
in the decomposition of $\d(a_{k+1})$ with respect to the basis $A$ with a nonzero coefficient.

\begin{lemma}
\label{dprimitive} Let $\d$ and $\d'$ be $(A,B)$-equivalent differentials. Then the set of all $\d$-trivial
elements coincides with the set of all $\d'$-trivial elements.
\end{lemma}

\pr Let $\d'=g\d g^{-1}$ for some $g\in \AutT(A,B)$. We denote by dots a linear combination of elements with
indices smaller than $k$.
We have $g^{-1}(a_{k+1})=\lambda_1a_{k+1}+...$, $g(a_{k})=\lambda_2 a_{k}+...$ for some $\lambda_1,\lambda_2
\ne 0$. Then $\d(a_{k+1})=\mu a_k+ ...$ with $\mu\ne 0$ implies $\d'(a_{k+1})=\lambda_2 \mu \lambda_1
a_k+...$.\qed

\subs{Set $\mathcal{D}_{A,B,C}$} We fix a triple $A\supset B\supset C$ of finite sets and suppose that $A$ is non-empty.
The set $A$ is graded and linearly ordered.
Denote by $\mathcal{D}_{A,B,G}$ the set of all $(A,B)$-differentials, such that any element of~$C$ is $\d$-trivial for any $\d \in \mathcal{D}_{A,B,C}$. 
For each (necesserily $\d$-trivial) element $a\in C$ we draw on the corresponding figure double segment joining $a$ with the next element in~$A$. 

Recall, that  for a strong Morse function $F$ we constructed (see Section~\ref{amodel}) a triple of graded sets $A_F\supset B_F\supset C_F$ ($C_F$ consists of all critical outward critical points).

\begin{statement}
\label{Morse_in_D} Differential of an algebraic model of a strong Morse function~$F$ belongs to~$\mathcal{D}_{A_F,B_F,C_F}$.
\end{statement}
\pr Consider an algebraic model $\M_{A_F,B_F,\d}$ of $F$.
%
We need to show that any element  $a_k\in C_F$ is  $\d$-trivial. By construction $a_{k+1} \in A_F\setminus B_F$.
Elements $a_k,a_{k+1}$ are generators of a complex calculating relative homologies of a pair $(F_{c+\varepsilon}, F_{c-\varepsilon})$, where $c=F(a_k)$. 
Since $F_{c-\varepsilon}$ is a strong deformation retract of $F_{c-\varepsilon}$ by Proposition~\ref{bifurcation}(3) these homologies are trivial and hence $a_k$ is $\d$-trivial. \qed

The space $\mathcal{D}_{A,B,C}$ is a main object of the paper. It is obviously a union of classes of equivalence. 
One can show that in general the number of classes of equivalence containing in $\mathcal{D}_{A,B,C}$ is infinite.

\subs{Quotient $M$-complex} Consider an $M$-pair $\M_{A,B,\,\d}$.
We will identify the set $A\setminus B$ with a basis of the quotient complex $\M_{A,\d}/\M_{B,\d_\lo{B}}$.
The linear order and grading on $A$ induce in a natural way the linear order and grading on $A\setminus B$.
It is clear that the induced differential on $\M_{A,\d}/\M_{B,\d_\lo{B}}$ is an $M$-differential with respect
to the linear order and grading on $A\setminus B$.
We denote the induced differential by $\d_{A\setminus B}$ and the quotient $M$-complex by $\M_{A\setminus
B,\d_\lo{A\setminus B}}$.

\subs{Weak equivalence}  We say that differentials $\d_1, \d_2\in \mathcal{D}_{A,B,C}$ are weakly equivalent if there exists an automorphism $g\in \Aut(A,B)$ ($g$ is not necessarily upper-triangular), such that $g\d_1=\d_2g$ and naturally induced automorphism  $g_\lo{A\setminus B}$ belongs to $\AutT(A\setminus B)$ and the restriction $g|_{\E\otimes B}$ belongs to $\AutT(B)$.

We will show that $\mathcal{D}_{A,B,C}$ is a disjoint union of weak equivalence classes. 
In each class of weak equivalence we will construct a unique convenient representative (a section of weak equivalence). 
To do that we need the following definitions and construction.

\subs{Direct sum decomposition} We say that an $M$-pair $\M_{A,B,\,\d}$ is {\it decomposable into a direct
sum} of two $M$-pairs (or, equivalently, a differential $\d$ is decomposable into direct sum), if there
exists a decomposition $A=A_1\cup A_2$ into disjoint nonempty subsets, such that the subspaces
$\E(A_1)$ and $\E(A_2)$ are $\d$-invariant.

In this case spaces $\E(A_i\cap B)$ are also $\d$-invariant, so the $M$-pairs $\M_{A_1,B_1,\,\d_1}$
and $\M_{A_2,B_2,\,\d_2}$ are well defined, where $B_i=A_i\cap B, i\in\{1,2\}$ and the $\d_i$ are the
restrictions of $\d$. 
We will write $\M_{A,B,\,\d}=\M_{A_1,B_1,\,\d_1}\oplus \M_{A_2,B_2,\,\d_2}$ in this
case. 
If $\M_{A,B,\,\d}=\M_{A_1,B_1,\,\d_1}\oplus \M_{A_2,B_2,\,\d_2}$ and $C$ is a set of $\d$-trivial
elements then, the set $C_i=B_i\cap C$ ($i\in\{1,2\}$) consists of $\d_i$-trivial elements and we write
$\M_{A,B,C,\,\d}=\M_{A_1,B_1,C_1,\,\d_1}\oplus \M_{A_2,B_2,C_2,\,\d_2}$ (or $\d=\d_1\oplus\d_2$) in that
case. 
Decomposition into the direct sum of greater number of summands is defined in a similar way. 
Decomposition of an $M$-pair into a direct sum of indecomposable summands is unique up to a reordering of the
summands. 

\ign{
The following lemma is obvious.

\begin{lemma}
Let $\M_{A,B,\,\d}=\M_{A_1,B_1,\,\d_1}\oplus ... \oplus \M_{A_K,B_K,\,\d_K}$. Then, $H(\d)$ =
$H(\d_1)\cup...\cup H(\d_K)$ and any $\d$-trivial element $a\in A$ is $\d_i$-trivial for some $i$. \qed
\end{lemma}
}

\subs{Operation $\#$} Consider the following partially defined operation on $M$-pairs. 
Consider non-zero $M$-pairs $\M_{A,B,C,\,\d}$ and $\M_{X,Y,Z,\,\delta}$, such that $A\cap X=\varnothing$. Let $A=\{a_1\prec...\prec a_K\}, X=\{x_1\prec...\prec x_L\}$. 
We suppose that degree of $x_L$ is bigger by one than degree of $a_1$ and
that $a_1\in B$, $a_1\notin C$, $x_L\in X\setminus Y$, $x_{L-1}\notin Z$. 
We denote by $A\#X$ the set $A\cup X$ with the order 
$$
x_1\prec...\prec
x_{L-1}\prec a_1\prec x_L \prec a_2 \prec ...\prec a_K.
$$  
We define a linear operator $\d\#\delta$ on elements of $A\#X$ as follows: $\d\#\delta(a_i)=\d(a_i)$ for any $i\in\{1,...,K\}$, $\d\#\delta(x_i)=\delta(x_i)$ for any $i\in\{1,...,L-1\}$ and $\d\#\delta(x_L)=\delta(x_L)+a_1$.
%
Denote by $C\#Z$ the set $C\cup Z\cup \{a_1\}$. 
Clearly, $\d\#\delta$ is an $M$-differential and the set $C\#Z$ consists of $\d\#\delta$-trivial elements. 
We denote $M$-pair $\M_{A\#X,B\#Y,C\#Z,\,\d\#\delta}$ ($B\#Y=B\cup Y$) as $\M_{A,B,C,\,\d}\#\M_{X,Y,Z,\,\delta}$. We set $\M\#0=\M$, $0\#\M=\M$ for any $M$-pair $\M$.
%
Operation $\#$ is obviously associative.

\begin{figure}[htb]
\begin{center}
\includegraphics[width=5 in]{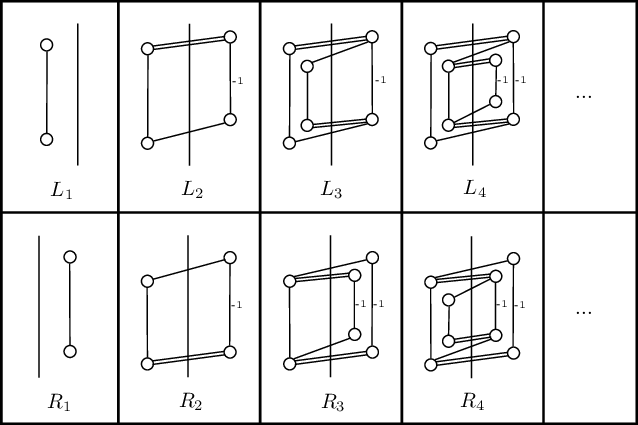}
\caption{\,\,Definition of $M$-pairs $L_k$ and $R_l$}
\label{Fig:RLdef}
\end{center}
\end{figure}

\subs{Definition of $M$-pairs $L_k$ and $R_l$} Let us define $M$-pairs $L_k$ and $R_l$ ($k,l\ge 0$) as follows: $L_0=R_0=0$, and for $k,l\ge 1$ we define $L_k, R_l$ (up to the common shift of grading) from the (infinite) table on the Fig.~\ref{Fig:RLdef}.  
The grading of any $M$-pair $L_k,R_l$ is defined uniquely up to a common shift. Each $L_k, R_l$ belongs to the corresponding $\mathcal{D}_{A,B,C}$ where $C$ is determined by double segments.

\subs{Main Theorem} Consider a triple $(A,B,C)$ of finite graded sets. We assume that $A$ is nonempty.
\begin{theorem} The number of classes of weak equivalence is finite. Any $M$-differential $\d\in \mathcal{D}_{A,B,C}$
 weakly equivalent to  unique differential decomposable into a direct sum of differentials of $M$-pairs of type:
\[
L_k\# R_l (k+l\ge 1), \,\,L_k\#\raisebox{-1.5ex}{\m{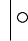}} (k\ge 0), \,\, \raisebox{-1.5ex}{\m{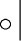}}\,\,\#R_l (l\ge 0),
\]
\[
L_k\#\raisebox{-1.5ex}{\m{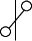}}\#R_l, \,\,\raisebox{-1.5ex}{\m{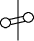}}, \,\,\raisebox{-1.5ex}{\m{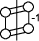}}.
\]
\label{th1}
\end{theorem}

\ign{

\subs{Totally decomposable $M$-differentials} Consider a picture corresponding to a $M$-pair with a
$M$-differential $\d \in \mathcal{D}_{A,B,G}$. Let element $a=a_k$ belong to the set $G$.
To emphasize this we will draw a double segment between the circle corresponding to $a$ and the circle
corresponding to the element $a_+=a_{k+1}$.
Elements of the set $H(\d)$ are shown as black circles.
The grading of any $M$-pair $\M_1$--$\M_{16}$ shown in Figure \ref{Fig:table} is defined uniquely up to a
common shift.
\ign{
\begin{figure}[htb]
\begin{center}
\includegraphics[width=5.5in]{vse_m.eps}
\caption{\,\,Glyphes}\label{Fig:table}
\end{center}
\end{figure}
}

We will call $M$-pairs of Figure~\ref{Fig:table} {\it glyphes}. We say that an $M$-differential $\d \in
\mathcal{D}_{A,B,G}$ is {\it totally decomposable} if $\d$ is a direct sum of glyphes.

\subs{The $M$-model Theorem}

\begin{theorem}
\label{main}{\rm(\,$M$-model Theorem)} For any $\d \in \mathcal{D}_{A,B,G}$ there exists a totally
decomposable $M$-differential {$G$-equivalent} to $\d$.

Moreover, for any $M_{A,B}$-differential $\d$ such $G$-equivalent differential may be chosen canonically.
Namely, for each triple $(A,B,G)$ exists a well defined map
$\mathcal{P}_{A,B,G}\colon\mathcal{D}_{A,B,G}\to\mathcal{D}_{A,B,G}$ satisfying the following properties:
\begin{enumerate}
\item $\mathcal{P}_{A,B,G}(\d)$ is totally decomposable;
\item for any $\d \in \mathcal{D}_{A,B,G}$, the differential $\mathcal{P}_{A,B,G}(\d)$ is $G$-equivalent to
$\d$;
\item if differentials $\d, \d' \in \mathcal{D}_{A,B,G}$ are weakly equivalent than
$\mathcal{P}_{A,B,G}(\d)=\mathcal{P}(\d')$;
\item if $\d$ is totally decomposable then $\mathcal{P}_{A,B,G}(\d)=\d$;
\item if $\M_{A,B,G,\,\d}=\M_{A_1,B_1,G,\,\d_1}\oplus \M_{A_2,B_2,G,\,\d_2}$ then
$\mathcal{P}_{A,B,G}(\d)=\mathcal{P}_{A_1,B_1,G_1}(\d_1)\oplus\mathcal{P}_{A_2,B_2,G_2}(\d_2)$.
%
\end{enumerate}
We say that $\mathcal{P}(\d)=\mathcal{P}_{A,B,G}(\d)$ is the {\it $M$-model} of the differential $\d$. The
map $\mathcal{P}$ is defined in the proof of the theorem.
\end{theorem}

}

\subs{Remark} The notion of weak equivalence, $G$-equivalence and last theorem are a product of naive
attempts of ``simplifying'' an $M$-differential. For closed manifold Theorem~

\section{On $M$-complexes and $M$-pairs}

This section contains preliminary results needed for the proof of Theorem~\ref{th1}. 
For a given $M$-differential we construct an equivalent $M$-differential having relatively simple form. 
We call differentials of this type quasi-elementary differentials. 
This construction is, in fact, the first step in the proof of Theorem~\ref{th1}. 
The proof of Theorem \ref{th1} (Sec.~\ref{Proof} below) involves subsequent simplifications of the
constructed quasi-elementary differential.
%
On the level of matrix elements of differential we choose an equivalent differential having many zeroes in its matrix.

\subs{On the structure of an $M$-complex}\label{Mcomstruct}
%
%
Let $A$  be a finite linearly ordered graded set. An $M_A$-differential $\partial$ is called an {\it
elementary differential}, if it satisfies  the following two conditions:
\begin{enumerate}
\item{for each $a\in A$ either $\partial(a)=0$, or there exists $b\in A$ such that $\partial(a)=b$;}
\item{$\partial(x)=\partial(y)=z$ and $x,y,z \in A$ implies $x=y$.}
\end{enumerate}
\ign{ Let $A=\{a_1\prec ...\prec a_K\}$  be a finite linearly ordered graded set. A $M_A$-differential
$\partial$ is called an {\it elementary $M$-differential}, if it satisfies  the following two conditions:
\begin{enumerate}
\item{for each $j\in\{1,\dots,\?K\}$ either $\partial(a_j)=0$, or there exists an index $m\in\{1,\dots,\?K\}$,
such that $\partial(a_j)=a_m$;}
\item{$\partial(a_j)=\partial(a_{j'})=a_m$ implies $j=j'$.}
\end{enumerate}
}
\begin{theorem}{\rm\cite{Serega}}     \label{Normal}
%
Any $M$-differential $\partial$ over a field $\E$ is equivalent to a unique elementary $M$-differential.
\hfill$\sq$
\end{theorem}

\pr An elementary $M$-differential, equivalent to the given one, is constructed explicitly by induction as
follows. For an one-dimensional $M$-complex the statement is trivial.
Assume that the statement is proved for $M$-complexes of dimension $k$. Let $A=\{a_1\prec ... \prec a_{k+1}\}$.
Consider an $M$-complex $\M_{A,\d}$. A
restriction $\d_1$ of the $M$-differential $\d$ on subcomplex $\mathcal{F}_k(A)$ could be reduced to an
elementary form by the induction hypothesis. We assume that this is done already.

We define a map $\d^{-1}$ on elements $a_1,...,a_k$ as follows: if $\d_1(a_j)=a_i$, then $\d^{-1}_1(a_i) =
a_j$, in other case $\d^{-1}_1(a_i) = 0$. Consider $\d(a_{k+1})$. If $\d(a_{k+1})=0$, then $\d$ is elementary
and the statement is proved. Let $\d(a_{k+1})=\sum\limits_{i\in I\subset \{1,\ldots,K\}}\lambda_i a_i$,
$\lambda_i \ne 0$ for $i \in I$. Since differential $\d_1$ is elementary differential and $\d^2(a_{k+1})=0$,
then for each $i \in I$ holds $\d(a_i)=0$. We decompose the set $I$ as a disjoint union of a set $I_1$,
consist from all $i\in I$, such that $\d^{-1}_1(a_i) \ne 0$ and a set $I_2$, consists of all other elements
of the set $I$.

Suppose that the set $I_2$ is empty. Consider an automorphism $T$, defined on basis elements as follows:
$a_{k+1}\mapsto a_{k+1}-\sum\limits_{i \in I_1} \lambda_i \d^{-1}_1(a_i))$, all rest basis elements are fixed
by $T$. Obviously $T\in \AutT(A)$. The differential $\d_T=T^{-1}\d T$ is an elementary $M$-differential,
since $\d_T(a_i)=\d(a_i)$ for $i\in\{1,\ldots,k\}$, and $\d_T(a_{k+1})=0$.

Suppose that the set $I_2$ is nonempty. 
Let $l$ be a maximal element of the set $I_2$. 
Consider an automorphism $T$, defined on basis elements: $a_{k+1}\mapsto a_{k+1}-\sum\limits_{i \in I_1}
\lambda_i\d^{-1}_1(a_i)$, $a_{l}\mapsto\sum\limits_{i\in I_2}\lambda_ia_i$, all rest basis elements are fixed
by $T$. 
Obviously $T\in \AutT(A)$. $M$-differential $\d_T$ is an elementary $M$-differential, since
$\d_T(a_{k+1})=a_m$, $\d^{-1}_1(a_m)=0$ and on all rest basis elements $\d_T$ coincide with $\d$.

Let us prove uniqueness. Denote by $d(m,n,\d)$ the dimension of relative homologies (see \cite{Ch-P}) 
$H_*(\mathcal{F}_m(A),\mathcal{F}_n(A),\d)$. Obviously $d(m,n,\d)$ depends only on a class of equivalence of
differential $\d$. Let $\d$ be an elementary $M$-differential $\d$ and  $\d(a_i)=a_j$. 
Then, obviously, the following equalities hold:
$$
d(i,j,\d)      = d(i-1,j-1,\d)      =   d(i-1,j,\d) \?+\? 1  = d(i,j-1,\d)  \?+\? 1.
$$
And conversely -- it is easy to check that if these equalities hold, than $\d(a_i)=~a_j$. Uniqueness is proved.
\qed

\subs{Partition of a basis of an $M$-complex into pairs and homologically essential elements}
Consider an $M$-complex $\M_{A,\d}$, $A=\{a_1\prec ... \prec a_N\}$.
Let $\d_1$ be the elementary $M$-differential equivalent to $\d$ (as in Theorem~\ref{Normal}).
We say that basis elements $a_i, a_j$ of an $M$-complex $\M_{A,\d}$ form a {\it $\d$-pair} if
$\d_1(a_i)=a_j$.
We say that a basis element is {\it ($\d$-)homologically essential} element if it does not appear in a
$\d$-pair.
According to Theorem~\ref{Normal}, any element of $A$ is either a homologically essential element or a member
of a unique $\d$-pair. 
%
The following assertions describe this combinatorial structure on $A$ in terms of an $M$-differential $\d$.

For $A=\{a_1\prec ... \prec a_N\}$ denote by $A^k$ the set $\{a_1\prec ... \prec a_k\}$. 
We denote by $\iota$ the inclusion $\E(A^j)\to \E(A)$ (the value of $j$ will be clear from the
context) and by $\iota_*$ the induced map in the homology.

\begin{lemma} \label{Charact}
(1) An element $a_j\in A$ of degree $l$ is homologically essential, if and only if
$$
\dim \iota_*(H_l(\E(A^j),\d))=\dim \iota_*(H_l(\E(A^{j-1}),\d))+1.
$$
The number of homologically essential elements of degree $k$ is equal to $\dim H_k(\M_{A,\d})$. \\
(2) An element $a_i\in A$ is not homologically essential if and only if 
$$
\iota_*H_*(\E(A^i),\d)=
\iota_*H_*(\E(A^{i-1}),\d).
$$
(3) Elements $a_m,a_n\in A$, $m>n$ form a $\d$-pair if and only if
\begin{equation*}
\begin{split}
&\dim H_*(\E(A^m),\E(A^n),\d)=\dim H_*(\E(A^{m-1}),\E(A^{n-1}),\d)= \\
& =\dim H_*(\E(A^{m-1}),\E(A^k),\d)+1=\dim H_*(\E(A^m),\E(A^{n-1}),\d)+1.
\end{split}
\end{equation*}
\end{lemma}

\pr The statement of the Lemma is obvious for an elementary differential.
%
Hence, by Theorem~\ref{Normal} it
holds for any $M$-differential, since the dimensions of homologies in the statement depend only on the
equivalence class of the $M$-differential. \qed

\subs{Boundary homologically essential elements of $M$-pairs}\label{grboundess} 
A filtration on a topological space or on a chain complex gives rise to a filtration on any subspace of its
homology.
This simple observation leads us to the following definition. 


Consider an $M$-pair $\M_{A,B,\,\d}$. Let $B=\{b_1\prec...\prec b_K\}$.
The subspace $\E(B^k) = \E\otimes\{b_1,...,b_k\}\subset \E(B)$ is $\d_\lo{B}$-invariant,
and hence graded space of homologies $H_*(\E(B^k),\d_\lo{B})$ is well defined.
Let $\iota_*\colon H_*(\E(B^k),\d_\lo{B})\to H_*(\E(B),\d_\lo{B})$ denote the map induced
by the inclusion $\E(B^k)\hookrightarrow \E(B)$.
Let $\d_*\colon H_{*+1}(\E(A),\E(B), \d) \to H_{*}(\E(B), \d_\lo{B})$ be the boundary map of the long
exact sequence of the pair $(\E(A),\E(B),\,\d)$.
%
Introduce the intersections
$$
I_k = \iota_*H_*(\E(B^k),\d_\lo{B})\cap \d_*H_*(\E(A),\E(B), \d) \subset
H_*(\E(B),\d_\lo{B}).
$$

A basis element $b_k\in B$ is called {\it ($\d$-)boundary homologically essential}, if $I_k \ne I_{k-1}$ (we
set $I_0 = 0$).

We will denote by $H(\d)$ the set of all $\d$-boundary essential elements. 
It is clear that if an $M$-differential $\d$ is $(A,B)$-equivalent to an $M$-differential $\d'$ then $H(\d)=H(\d')$.

\subs{Homologically essential and boundary homologically essential basis elements} Consider an $M$-pair
$\M_{A,B,\d}$.
Some elements of the set $B$ are $\d$-boundary homologically essential (see Section~\ref{grboundess}).

\begin{lemma}
\label{boundessent} Every $\d$-boundary homologically essential element is \\$\d_{B}$-homologically essential.
The number of $\d$-boundary homologically essential elements of degree $k$ is equal to
$\dim\d_*(H_{k+1}(\E(A),\E(B),\d))$.
\end{lemma}

\pr According to Lemma~\ref{Charact} applied to the $M$-complex $\M_{B,\d_{B}}$,
$\iota_*H_*(\E(B^l),\d_\lo{B})\ne \iota_*H_*(\E(B^{l-1}),\d_\lo{B})$ if and only if $b_l$ is $\d_\lo{B}$-homologically essential.
For such $l$ the dimensions of the spaces $\iota_*H_*(\E(B^l),\d_\lo{B})$ and
$\iota_*H_*(\E(B^{l-1}),\d_\lo{B})$ differ by $1$, so we get a full graded flag in
$\iota_*H_*(\E(B),\d_\lo{B})$, consisting of spaces $\iota_*H_*(\E(B^k),\d_\lo{B})$.
The intersection of a graded subspace with a full flag is a full flag in the subspace.
This proves the first claim of Lemma.
Spaces $\iota_*H_*(\E(B^l),\d_\lo{B})$ and $\d_*(H_*(\E(A),\E(B),\d))$ are direct sums of their homogeneous components. This proves the second claim.\qed

\ign{
\subs{Quotient $M$-complexes} Consider an $M$-pair $\M_{A,B,\,\d}$.
We will identify the set $A\setminus B$ with a basis of the quotient complex $\M_{A,\d}/\M_{B,\d_\lo{B}}$.
The linear order and grading induce in a natural way the linear order and grading on $A\setminus B$.
It is clear that the induced differential on $\M_{A,\d}/\M_{B,\d_\lo{B}}$ is an $M$-differential with respect
to the linear order and grading on $A\setminus B$.
We denote the induced differential by $\d_{A\setminus B}$ and the quotient $M$-complex by $\M_{A\setminus
B,\d_\lo{A\setminus B}}$.
An isomorphism of $M$-pairs naturally induces an isomorphism of quotient $M$-complexes.
}

\subs{The sets $P,Q,R$ and $X,Y,Z$ and bijections}\label{thesets}
%
Consider an $M$-pair $\M_{A,B,\,\d_0}$. By Theorem~\ref{Normal}, the $M$-differential $\d_0$ of the $M$-pair
$\M_{A,B,\,\d_0}$ is $(A,B)$-equivalent to an $M$-differential $\d$, such that $\d_\lo{B}$ and
$\d_\lo{A\setminus B}$ are elementary $M$-differentials.

The differential $\d_\lo{B}$ is elementary, hence, the set $B$ is decomposed into a disjoint union of subsets
$P,Q,R$ such that $\d_\lo{B}$ restricts to a bijection $Q\to R$ and $\d_\lo{B}$ to zero on $P$.
Similarly, the differential $\d_\lo{A\setminus B}$ is elementary, hence, the set $A\setminus B$ is decomposed
into a disjoint union of subsets $X,Y,Z$ such that $\d_\lo{A\setminus B}$ restricts to a bijection $Y\to Z$
and $\d_\lo{A\setminus B}$ to zero on $X$.

By Theorem~\ref{Normal}, the sets $P,Q,R$ and $X,Y,Z$ and bijections $Q\to R$ and $Y\to Z$ depend only on the
equivalence class of the differential $\d_0$.

\subs{The definition of a quasi-elementary differential}\label{quasiel} Consider a vector space $L$ and a
basis $W$ of $L$.
We say that $v=\sum_{w\in W}v_ww\in L$ contains an $a\in W$ (or $a$ appears in $v$) if $v_a\ne 0$.

We say that an $M_{A,B}$-differential $\d$ is {\it quasi-elementary} if it satisfies following conditions:
\begin{enumerate}
\item{differentials $\d_\lo{B}$ and $\d_\lo{A\setminus B}$ are elementary;}
\item{for each element $x\in X$, the vector $\d(x)$ contains at most one element of $P$;}
\item{if $\d(x)$ contains $p\in P$ then the corresponding coefficient $\d(x)_p$ equals to 1;}
%
\item{any element of the set $P$ appears in at most one vector $\d(x)$ for $x\in X$.}
\end{enumerate}

\subs{From $M$-differentials to quasi-elementary differentials}\label{quasiel1} Recall that $H(\d)$ denotes
the set of $\d$-boundary homologically essential elements of a $M$-differential $\d$ (see
Section~\ref{grboundess}).
\begin{lemma}
\label{quasisimple} $(1)$ Any $M$-differential is equivalent to a quasi-elementary $M$-differential.\\
$(2)$ Suppose that $\d$ is a quasi-elementary $M$-differential.  The set $H(\d)$ coincides with the set of
all elements of $P$ appearing in vectors $\d(x)$ for $x \in X$.\\
$(3)$ Suppose that $\d$, $\d_1$ are equivalent quasi-elementary differentials. For any $x\in X$ and $p\in P$,
$\d(x)$ contains $p$ if and only if  $\d_1(x)$ contains $p$.
\end{lemma}

Note that an equivalence class of $M$-differentials may contain more than one quasi-elementary differential.
The proof of Lemma \ref{quasisimple} is given in~\ref{proofqs} below.

\begin{wrapfigure}{r}{75pt}
\begin{center}
\includegraphics[width=.2\textwidth]{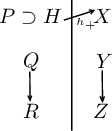}
 \end{center}
 \caption{}
\label{PQR}
\end{wrapfigure}
\subs{The injection $h_+$} Consider an $M$-pair $\M_{A,B,\d}$. 
We define the map $h_+\colon H(\d) \to X$ as
follows. 
Let $\d'$ be a quasi-elementary differential equivalent to $\d$.
For  $b\in H(\d)$, $h_+(b)$ is, by definition, the element $x\in X$ such that $\d'(x)$ contains $b$.
By Lemma~\ref{quasisimple}, the map $h_+$ is defined correctly and depends only on the equivalence class of
$\d$. It is clear that $h_+$ is an injection. Notice  that $h_+$ increases degree by~$1$.

The following assertion is an immediate corollary of Lemma~\ref{quasisimple}.
\begin{corollary}
\label{PQRDEF} 
The partitions
$$
B=P\sqcup Q\sqcup R, A\setminus B = X\sqcup Y\sqcup Z,
$$
the bijections $Q\to R$ and $Y\to Z$, the subset $H=H(\d)\subset P$ and the injection $h_+\colon H \to X$ are
invariants of an $M$-pair $\M_{A,B,\d}$ (see Fig.\ref{PQR} for schematic diagram).\qed
\end{corollary}

\subs{Proof of Lemma \ref{quasisimple}} \label{proofqs} We prove the first claim of the Lemma. 
Consider an $M$-pair $\M_{A,B,\,\d}$, $A=\{a_1\prec ... \prec a_N\}$.

We assume that the differentials $\d_\lo{B}$, $\d_\lo{A\setminus B}$ are elementary.
The claim is obvious if the set $X$ is empty. Suppose now that $X$ is not empty. Denote its elements
$x_1<...<x_I$ with respect to the order induced from the order on $A$. 
We will use induction to prove the
following statement: $\d$ is equivalent to a differential $\delta=\delta_k$ such that
\begin{enumerate}
\item{for any $i\in \{1,...,k\}$ vector $\delta(x_i)$ contains at most one element from $P$;}
\item{each element from $P$ appears in at most one vector $\delta(x_j)$ for $j \in \{1,...,k\}$ and
$\d_\lo{B}=\delta_\lo{B}$, $\d_\lo{A\setminus B}=\delta_\lo{A\setminus B}$.}
\end{enumerate}
\ign{ for any $i\in \{1,...,k\}$ vector $\delta(d_i)$ contains at most one element from $P$, each element
from $P$ appears in at most one vector $\delta(d_j)$ for $j \in \{1,...,k\}$ and $\d_\lo{B}=\delta_\lo{B}$,
$\d_\lo{A\setminus B}=\delta_\lo{A\setminus B}$. }

Let $k=1$. The vector $\d(x_1)$  may contain only elements from sets $P,R$ since $\d_\lo{A\setminus
B},\d_\lo{B}$ are elementary and $\d^2(x_1)=0$. 
Hence, $\d(x_1)=p+r$, $p\in \E(P), r\in\E(R)$. If $p=0$ then the first step of the induction is proved. 
Consider the case $p\ne 0$. Let
$p_i$ be the maximal element of $P$ which appears in $p$. 
Consider $T_p\in \AutT(A,B)$ such that $T_p(p_i)=p$ and the rest of basis elements is fixed by $T_p$. 
Then, the differential $\delta_1=T_p^{-1}\d T_p$ has the properties desired. 
This establishes the base of the induction.

Suppose that $\d$ is equivalent to $\delta_k$ as above. We may assume that $\d=\delta_k$.
If $\d(x_{k+1})$ does not contain elements from $P$ then $\delta_{k+1}=\d$ is the differential we need.
Let $x \in \E(\{x_1,...,x_{k}\})$. Denote by $T_x\in \AutT(A,B)$ the automorphism which maps
$x_{k+1}$ to the sum $x_{k+1}+x$ and fixes all other basis elements.
For every $x \in \E(\{x_1,...,x_{k}\})$, the differential $T_x^{-1}\d T_x$ satisfies the $k$-th
induction hypothesis and, for a suitable $x_0 \in \E(\{x_1,...,x_{k}\})$, the vector $T_{x_0}^{-1}\d
T_{x_0}(x_{k+1})$ does not contain elements from $P$ which appear in $T_{x_0}^{-1}\d T_{x_0}(d_i)$ for
$i\in\{1,...,k\}$. 
So, we may assume that $\d(x_{k+1})$ does not contain elements from $P$ which appear in
$\d(x_i)$ for $i\in\{1,...,k\}$. Let $\d(x_{k+1})=p+q$. 
Then, the differential $\delta_{k+1}=T_p^{-1}\d T_p$ has the properties desired. The first claim of Lemma is proved.

Let $\d$ be a quasi-elementary differential. 
Then, space $I_k= \iota_*H_*(\E(B_k),\d_\lo{B})\cap \d_*H_*(\E(A),\E(B), \d)$ from the definition
of $\d$-boundary homologically essential elements is spanned on elements
$b_i\in\{b_1,...,b_k\}\cap P$ appearing in vectors $\d(x)$ for $x\in X$. 
%
%
This proves the second claim of Lemma.


We are now going to prove the third statement. Let $(A,B)$ be a pair of finite ordered graded sets. Let
$B=\{b_1\prec...\prec b_K\}$, $A\setminus B =\{a_1\prec...\prec a_L\}$.  Denote by $A'$ the ordered graded
set which is equal (as a graded set) to $A$ with the linear order $\prec_n$:
$$
b_1\prec_n...\prec_n b_K\prec_n a_1\prec_n ...\prec_n a_L.
$$
An $M_{(A,B)}$-differential $\d$ naturally induces the $M_{A'}$-differential $\d'$ which is equal (as a linear map) to $\d$. 
Obviously, the map sending $\d$ to $\d'$ maps $(A,B)$-equivalent differentials to $A'$-equivalent ones.
%

In view of Theorem~\ref{Normal}, the third statement follows from the following assertion.

\begin{lemma} Let $\d$ be a quasi-elementary differential. If $p\in H(\d)$ appears in $\d(x)$ for $x\in X$
then the elements $x,p$ form a $\d'$-pair.
\end{lemma}

\pr From the assumption of the Lemma it follows that $\d(x)=\lambda p+r$, where $\lambda\ne 0$ and $r\in
\E(R)$. Let $q\in\E(Q)$ satisfy $\d(q)=r$. 
It easily follows from $\d^2=0$ that the element $p\in H(\d)$ does not appear in any $\d(z), z\in Z$.
Let $y_1,...,y_k\in Y$ be all elements such that $p$ appears in $\d(y_i)$ with a coefficient $\lambda_i\ne 0$. 
Denote $\d_\lo{A\setminus B}(e_i)$ by $z_i$.
Consider a map $T$ such that $T(z_i)=z_i+\lambda_i p$ for $i\in\{1,...,k\}$ and $T$ fixes all other basis
elements. Each $z_i$ is bigger than $p$ in the set $A'$.
Hence, $T \in \AutT (A')$. The differential $\d'_1=T^{-1}\d' T$ is $M_{A'}$-differential $A'$-equivalent to $\d'$.
The element $p$ appears in $\d'_1(x)$ and does not appear in the $\d'_1$-images of other basis elements.
Denote by $T_1$ an automorphism, such that $T_1(x)=x-q$ and $T_1$ fixes all other elements.
The differential $\d'_2=T_1^{-1}\d'_1T_1$ is $A'$-equivalent to $\d'_1$ and the subspace $\E\otimes\{x,p\}$
is its direct summand. 
Now $\d'_2(x)=\lambda p$ means that $d,p$ form a $\d'$-pair.\qed

\section{The proof of Theorem~\ref{th1}} \label{Proof}

In this section we are proving Theorem~\ref{th1}.
Our strategy is the following. On the first step for a $M$-differential $\d \in \mathcal{D}_{A,B,C}$ we consider a quasi-elementary differential $\widetilde{\d} \in \mathcal{D}_{A,B,C}$. 
%
On the second step we start the procedure 
of eliminating of matrix elements of $\d'$. When the procedure stops we get a weakly equivalent to $\d$ quasi-elementary differential which is decomposable (up to the action of diagonal subgroup in $\AutT(A,B)$) into the summands of Theorem~\ref{th1}.


\subs{Nonvanishing matrix elements} Consider an $M$-differential
$\d\in\mathcal{D}_{A,B,C}$. According to Lemma~\ref{quasisimple}, there exists a quasi-elementary
differential $\widetilde{\d}$ equivalent to $\d$. By Lemma~\ref{dprimitive}, all elements of the set $C$ are
$\widetilde{\d}$-trivial. Thus $\widetilde{\d}\in\mathcal{D}_{A,B,C}$. We denote by $P,Q,R$, $X,Y,Z$, $H$ and
$h_+$ the sets and map from Statement~\ref{PQRDEF}. Let $\langle.,.\rangle$ be the standard scalar product on
$\E(A)$: $\langle a_i,a_j\rangle=\delta_{ij}$, the Kronecker symbol.

We will use the following simple Lemma. We leave its proof (which is easily following from $\widetilde{\d}^2=0$) to the reader.

\begin{lemma}\label{paral} Suppose $\widetilde{\d}$ is a quasi-elementary differential.

Assume elements $q\in Q$, $r\in R$, $y\in Y$, $z\in Z$ are such that $\widetilde{\d}(q)=r,
\widetilde{\d}_{A\setminus B}(y)=z$. Then $\langle \widetilde{\d}(y), q\rangle=-\langle\widetilde{\d}(z),
r\rangle$.

Assume elements $y\in Y$ and $z\in Z$ are such that $\widetilde{\d}_\lo{A\setminus B}(y)=z$, and $b\in B$
appears in $\widetilde{\d}(z)$. Then $b\in R$ and there exists $q \in Q$ such that $\widetilde{\d}(q)=b$ and
$q$ appears in $\widetilde{\d}(y)$. \qed
\end{lemma}

For an element $c\in C$ we denote by $c_+ \in A$ the minimal element bigger than $c$.
The element $c_+$ belongs to the set $A\setminus B$.
We also used notation $C_+$ for  the set of elements $c_+$ for all  $c\in C$.
%
We say that elements $c$ and $c_+$ form a $C$-pair (or $(c,c_+)$ is a $C$-pair). Also for element $a\in C_+$ we denote by $a_-$ an element $c$ such that $c_+=a$. 

Next Lemma describes a (distinguished) subset of non-zero coefficients of the matrix of $\widetilde{\d}$.

\begin{lemma}\label{nonzero}
Suppose that elements $a \in A\setminus B$, $b \in B$, $b\prec a$ satisfy at least one of the following
conditions:
\begin{enumerate}
\item{$(b,a)$ is a $C$-pair;}
\item{$b\in H$ and $a=h_+(b)$;}
\item{$a\in Y$, $b \in Q$, and elements $\widetilde{\d}(b)$ and $\widetilde{\d}{\lo{A\setminus B}}(a)$
generate a $C$-pair;}
\item{$a\in Z$, $b \in R$, and elements $q\in Q, y\in Y$ such that $\widetilde{\d}(q)=b$,
$\widetilde{\d}_{\lo{A\setminus B}}(y)=a$ generate a $C$-pair.}
\end{enumerate}
Then $\langle \widetilde{\d}(a),b\rangle \ne 0$.
\end{lemma}

\pr If $(a,b)$ satisfies to either (1) or (2) then $\langle \widetilde{\d}(a),b\rangle \ne 0$ by the
definitions of $C$ and $h_+$ respectively.
%
If $(a,b)$ satisfies (3) then, by Lemma~\ref{paral}, $\langle \widetilde{\d}(a),b\rangle=-\langle
\widetilde{\d}(\widetilde{\d}{\lo{A\setminus B}}(a)),\widetilde{\d}(b)\rangle$.
Since $(\widetilde{\d}(b),\widetilde{\d}{\lo{A\setminus B}}(a))$ is a $C$-pair we have $\langle \widetilde{\d}(a),b\rangle \ne 0$. 
For condition (4) we have $\langle \widetilde{\d}(y),q\rangle=-\langle \widetilde{\d}(a),b\rangle$ and, by Lemma~\ref{paral}, $\langle \widetilde{\d}(a),b\rangle \ne 0$ since $(q,y)$ is a $C$-pair.

\subs{Elimination of matrix elements. From an $M$-differential to a minimal differential} 
 We say that a quasi-elementary differential $\widetilde{\d}\in \mathcal{D}_{A,B,C}$ 
is {\it minimal}, if any pair $a \in A\setminus B$, $b \in B$, $b\prec a$ such that $\langle \widetilde{\d}(a),b\rangle \ne 0$ satisfies at least one of the conditions of Lemma~\ref{nonzero}. 
We call the corresponding $M$-pair minimal as well.

We call  differentials $\d,\d'$ {\it similar} if $\langle \d(a_i),a_j\rangle=0$ if and only if $\langle \d'(a_i),a_j\rangle=0$ for any $i,j$ (in other words the matrices of $\d$ and $\d'$ have zeroes at same places). 

\begin{lemma}\label{minimal} For any differential $\d \in \mathcal{D}_{A,B,C}$ there exists a minimal
differential $\delta\in \mathcal{D}_{A,B,C}$ weakly equivalent to $\d$. All such minimal differentials are
similar.

If differentials $\d,\d'\in \mathcal{D}_{A,B,C}$ are similar and $\d,\d'$ are weakly equivalent
(respectively) to minimal differentials $\delta,\delta'$ then $\delta,\delta'$ are similar.
\end{lemma}

\pr Let $\d_1$ be a quasi-elementary differential equivalent to $\d$. We will construct a finite sequence of
differentials $\d_1$,...,$\d_m$ such that for any $a\in A\setminus B, b\in B$ the equality $\langle
\d_{i-1}(a),b \rangle = 0$ implies $\langle \d_{i}(a),b \rangle = 0$ for any $i\in \{2,...,m\}$, differential
$\d_i$ is quasi-elementary and weakly equivalent to $\d_{i-1}$, and $\d_m=\delta$ is a minimal differential.
We proceed by induction. Suppose that $\d_i=\rho$ is not a minimal differential, consider a pair $a\in
A\setminus B$, $b\in B$ such that $\langle\rho(a),b\rangle\ne 0$ and $(a,b)$ does not satisfy to any of the
condition $(1)-(4)$ of Lemma~\ref{nonzero}.

\begin{lemma}\label{5cases} The pair $(a,b)$ satisfies to one of the conditions:

N1. $a\in Y$,  $b\in Q$ and each pair $(b,a)$, $(r,z)$ is not a $C$-pair, where $r=\rho(b),
z=\rho_\lo{A\setminus B}(a)$;

N2. $a\in Y$,  $b\in R$ and $(b,a)$ is not a $C$-pair;

N3. $a\in Y$, $b\in P$ and $(b,a)$ is not a $C$-pair;

N4. $a\in Z$,  $b\in R$ and  each pair $(b,a)$, $(q,y)$ is not a $C$-pair, where $q\in Q$ and $y\in Y$ be
such elements that $\rho(q)=b$, $\rho_\lo{A\setminus B}(y)=a$;

N5. $a\in X$, $b\in R$ and $(b,a)$ is not a $C$-pair.
\end{lemma}
\pr Consider the following cases $a\in Y$, $a\in Z$ and $a\in X$.

Let $a\in Y$. If $b\in Q$ then $(b,a)$ does not satisfy to conditions $(1)$, $(2)$ and $(4)$ of Lemma~\ref{nonzero} automatically.
Hence $N1$ holds, because it contradicts to the condition $(3)$. 
If $b\in R$ then $(b,a)$ does not satisfy to conditions $(2)$, $(3)$ and $(4)$ automatically. 
Hence $N2$ holds, because it contradicts to $(1)$. 
If $b \in P$ then $(y,x)$ does not satisfy to conditions $(2)$, $(3)$ and $(4)$ automatically and $N3$ holds.

Let $a\in Z$. Then, by Lemma~\ref{paral}, $b\in R$.  
The pair $(b,a)$ does
not satisfy to conditions $(2)$ and $(3)$ automatically. Hence it satisfies $N4$ since it contradicts to
conditions $(1)$ and $(4)$.

Let $a\in X$. In that case $b \notin P$ otherwise we have $a=h_+(b)$. The equality  $\rho^2(a)=0$ implies
$b\notin Q$. Hence $y\in R$ and we get the condition $N5$. \qed

For each case $N1$-$N5$ we define a subsequent differential $\d_{i+1}$ by the following Lemma.

\begin{lemma}\label{simplification} Let $\d_{i+1}$ be a linear operator such that:

If $(a,b)$ satisfies $N1$ then $\d_{i+1}(a)=\rho(a)-\langle\rho(a),b\rangle b$,
$\d_{i+1}(z)=\rho(z)-\langle\rho(z),r\rangle r$, where $r=\rho(b), z=\rho_\lo{A\setminus B}(a)$, and
$\d_{i+1}$ coincides with $\rho$ on all other basis elements;

If $(a,b)$ satisfies $N4$ then $\d_{i+1}(a)=\rho(a)-\langle\rho(a),b\rangle b$,
$\d_{i+1}(y)=\rho(y)-\langle\rho(y),q\rangle q$, where $q\in Q$ such that $b=\rho(q), a=\rho_\lo{A\setminus B}(y)$, and
$\d_{i+1}$ coincides with $\rho$ on all other basis elements;

If $(a,b)$ satisfies either $N2$ or $N3$ or $N5$ then $\d_{i+1}(a)=\rho(a)-\langle\rho(a),b\rangle b$ and
$\d_{i+1}$ coincides with $\rho$ on all other basis elements;

In each these cases $\d_{i+1}$ is quasi-elementary differential, $\d_{i+1} \in \mathcal{D}_{A,B,C}$ and
$\d_{i+1}$ is weakly equivalent to $\rho$.
\end{lemma}
\pr If $(a,b)$ satisfies $N1$ then $\d_{i+1}$ is equal to $S_1\rho S_1^{-1}$ for an automorphism $S_1$ such
that $S_1(x)=x-\langle\rho (a),b \rangle b$ and $S_1$ fixes all other basis elements. 
Thus $\d_{i+1}$ is differential ($\d_{i+1}^2=0$). 
Matrix elements $\langle\d_{i+1}(a),b\rangle$ and
$\langle\d_{i+1}(z),r\rangle$ are equal to zero. Indeed
$$
\langle\d_{i+1}(a),b\rangle=\langle\rho(a)-\langle\rho(a),b\rangle
b,b\rangle=\langle\rho(a),b\rangle-\langle\rho(a),b\rangle=0
$$
and
$$
\langle\d_{i+1}(z),r\rangle=\langle\rho(z+\langle\rho(a),b\rangle b)
,r\rangle=\langle\rho(a),b\rangle+\langle\rho(z),r\rangle=0.
$$
Obviously, all other matrix elements $\langle \d_{i+1}(a_m), a_n\rangle$ of $\d_{i+1}$ coincides with matrix
elements $\langle \rho(a_m), a_n\rangle$ of $\rho$. 
Hence the differential $\d_{i+1}$ belongs to the space $\mathcal{D}_{A,B,C}$ and it is quasi-elementary differential. 
The subspace $\E(B)$ is $S_1$-invariant and automorphisms $S_1|_B$ and $S_1|_{A\setminus B}$ are upper-triangular by construction, hence $\d_{i+1}$ is weakly equivalent to $\rho$.

Proof for other cases is analogous to the considered one.  If $(a,b)$ satisfies $Nj$ ($j\in\{2,3,4,5\}$) then
$\d_{i+1}=S_j \rho S^{-1}_j$. 
Operator $S_4$ acts nontrivially on the element $a$ only and $S_4(a)=a+\langle\rho(a),b\rangle q$.  
Operators $S_2$ and $S_5$ act nontrivially on $a$ only:
$S_2(a)=S_5(a)=a+\langle\rho(a),b\rangle q$, where $q\in Q$ and $\d(q)=b$. 
An operator $S_3$ acts nontrivially only on element $y=\rho_{A\setminus B}(a)$ and $S_3(y)=y-\langle\rho(a),b\rangle q$.
\qed

A set of non-zero matrix elements of minimal differential weakly equivalent to $\d$ is uniquely defined by
sets $P,Q,R,X,Y,Z,H$ and the map $h_+$. Hence all minimal differentials weakly equivalent to $\d$ are similar.

If differentials are weakly equivalent then sets $P,Q,R,X,Y,Z,H$ and the map $h_+$ coincide. It proves the
second claim of Lemma.\qed

\subs{End of proof. Indecomposable summands of minimal differential} The proof of Theorem~\ref{th1} ends with next Lemma.
\begin{lemma}\label{summands}
Assume that $\d\in \mathcal{D}_{A,B,C}$ is minimal quasi-elementary differential. 
Then any indecomposable summand of $\d$ is equivalent  to unique summand of Theorem~\ref{th1}.
\end{lemma}

\begin{figure}[htb]
\begin{center}
\includegraphics[width=5 in]{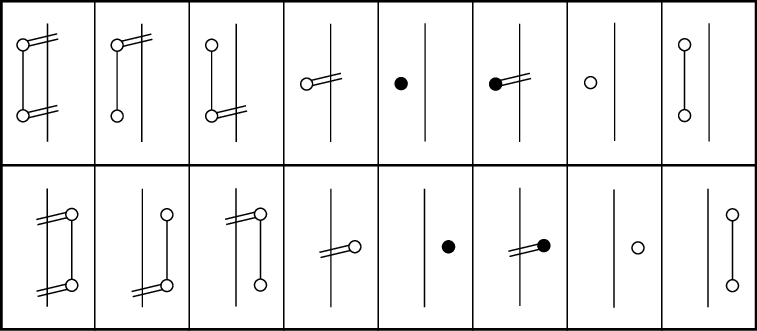}
\caption{\,\,Vertices of $\mathcal{G}(\d)$}
\label{Fig:VerticesofG}
\end{center}
\end{figure}

\pr  Let $\d\in \mathcal{D}_{A,B,C}$ be a quasi-elementary differential. 
Consider corresponding sets $P,Q,R,X,Y,Z,H$ and the map $h_+$. 
We define a graph $\mathcal{G}=\mathcal{G}(\d)$ as follows. 
Vertices of $\mathcal{G}$ are all one elements subsets in $P$ and in $X$, all two elements subsets $\{p,q\}$, such that $p \in P, q\in Q$ and $\d(p)=q$, and  all two element subsets $\{y,z\}$, such that $y \in Y, z\in Z$ satisfying 
$\d|_\lo{A\setminus B}(y)=z$.
We draw vertices corresponding to elements in $C$ and $C_+$
by circles with two segments, we draw vertices corresponding to elements of $H$ and $h_+(H)\subset X$ as black circles. 
It is easy to check that there are sixteen types of vertices.  All various possibilities are shown on Fig.~\ref{Fig:VerticesofG}. Vertices of $\mathcal{G}(\d)$ naturally correspond to parts in a picture for $\d$. 

Two distinct   vertices $\alpha$ and $\beta$ of $\mathcal{G}(\d)$ are connected by at most one edge and we connect $\alpha$ with $\beta$ by an edge if and only if there are elements $a\in \alpha$ and $b \in \beta$ such that  $\<\d(a),b\>\ne 0$ or  $\<\d(b),a\>\ne 0$. 

Obviously, connected components of $\mathcal{G}(\d)$ are in one to one natural correspondence with indecomposable summands of $\d$. By simple checking 
of minimality conditions we get the next Lemma.

\begin{lemma}
\label{isolatedvert}
Let $\d\in \mathcal{D}_{A,B,C}$ be a minimal quasi-elementary differential.
%
The corresponding equivalent direct summands of Theorem~\ref{th1} are $\raisebox{-1.0ex}{\m{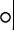}}\,\# R_0$, $L_1\# R_0$, $L_0\#\raisebox{-1.0ex}{\m{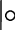}}$ and $L_0\# R_1$ respectively. \qed
\end{lemma}

In what follows we will assume that $\d\in \mathcal{D}_{A,B,C}$ is a minimal quasi-elementary differential and $\mathcal{G}=\mathcal{G}(\d)$ is its graph. 

Consider vertices of types $\raisebox{-1.0ex}{\m{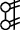}}$ and $\raisebox{-1.0ex}{\m{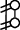}}$. 
We say that a vertex $\{p,q\}, p\in P, q\in Q, \d(p)=q$ of the type $\raisebox{-1.0ex}{\m{vgl1}}$ (resp. $\{y,z\}$, $y\in Y, z\in Z, \d(y)=z$ of the type $\raisebox{-1.0ex}{\m{vgr1}}$) is {\it closed} if elements $p_+,q_+$ form $\d_\lo{A\setminus B}$-pair  (resp. $y_-,q_-$ form $\d_B$-pair). 
Similarly, we say that a vertex of type \raisebox{-1.0ex}{\m{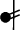}} (\raisebox{-1.0ex}{\m{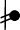}}) is {\it closed} if it consists of an element $c\in C\cap H$ such that $c_+=h_+(c)$ (resp. $c\in C_+\cap h_+(H)$ such that $h_+(c_-)=c$).
We say that a vertex of types \raisebox{-1.0ex}{\m{vgl1}}, \raisebox{-1.0ex}{\m{vgr1}}, \raisebox{-1.0ex}{\m{vgl6}} or \raisebox{-1.0ex}{\m{vgr6}} is {\it open} if it is not closed.

\begin{lemma}
\label{2,4}
Every closed vertex of the type \raisebox{-1.0ex}{\m{vgl1}} or  \raisebox{-1.0ex}{\m{vgr1}} corresponds to a direct summand of $\d$ which is equivalent to a summand  \raisebox{-1.5ex}{\m{s4}} of Theorem~\ref{th1}.
Every closed vertex of the type \raisebox{-1.0ex}{\m{vgl6}} or  \raisebox{-1.0ex}{\m{vgr6}} corresponds to a direct summand of $\d$ which is a summand  \raisebox{-1.5ex}{\m{s5}} of Theorem~\ref{th1}.
\qed
\end{lemma}

From now on we will assume that $\mathcal{G}=\mathcal{G}(\d)$ ($\d \in \mathcal{D}_{A,B,C}$) does not contain closed vertices and isolated vertices. 
For a vertex $a$ of $\mathcal{G}$ consider its closed neighbourhood $B(a)$ ($a$ and all vertices adjacent to $a$ in $\mathcal{G}$). 
We draw all elements of $A$ generating vertexes of $B(a)$ in a manner similar to drawing of $\d$ taking care on non-zero coefficients only. So the picture for $B(a)$ is a fragment for the picture of $\d$ up to values of non-zero coefficients on segments.

\begin{figure}[htb]
\begin{center}
\includegraphics[height=4.5 in, keepaspectratio]{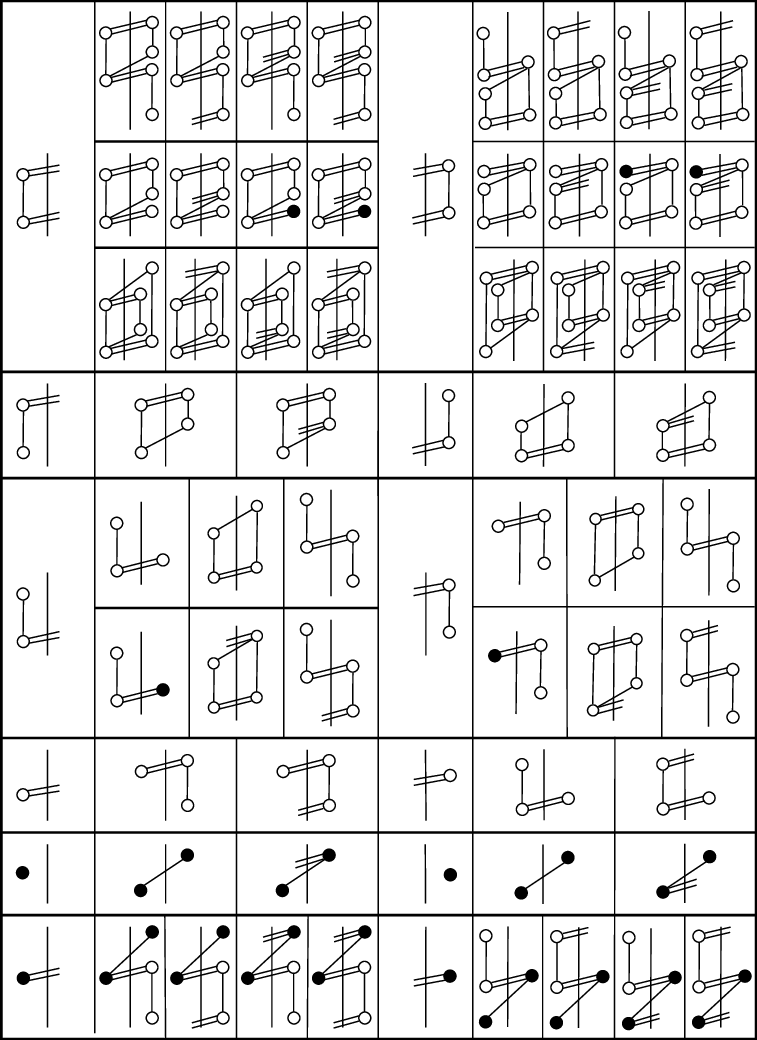}
\caption{\,\, Adjacent vertices of $\mathcal{G}$ in $\d$}
\label{Fig:Verticesedges}
\end{center}
\end{figure}

One can show by exhaustion that the table from  Fig.~\ref{Fig:Verticesedges} contains the full list of all possible pictures for $B(a)$. 
From that table we get the following statement.

\begin{statement} 
\label{val} Every vertex of $\mathcal{G}$ has valency 2 or 1. All vertices of type \raisebox{-1.0ex}{\m{vgl1}}, \raisebox{-1.0ex}{\m{vgr1}}, \raisebox{-1.0ex}{\m{vgl6}} and \raisebox{-1.0ex}{\m{vgr6}}  has valency $2$ in $\mathcal{G}$. 
All vertices of type \raisebox{-1.0ex}{\m{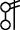}}, \raisebox{-1.0ex}{\m{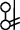}}, \raisebox{-1.0ex}{\m{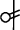}}, \raisebox{-1.0ex}{\m{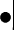}}, \raisebox{-1.0ex}{\m{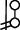}}, \raisebox{-1.0ex}{\m{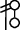}}, \raisebox{-1.0ex}{\m{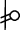}}, \raisebox{-1.0ex}{\m{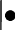}} have valency $1$ in $\mathcal{G}$.
\qed
\end{statement}

Accordingly to Fig.~\ref{Fig:Verticesedges} every vertex $v$ of type \raisebox{-1.0ex}{\m{vgl1}} or \raisebox{-1.0ex}{\m{vgr1}} has valency $2$ in $\mathcal{G}$. 
One of the edges beginning at $v$ is naturally distinguished and we equip it by orientation as follows.

The set $A$ is ordered, $A=\{a_1\prec ... \prec a_N\}$. Consider an auxiliary order $\prec_C$ on $A$ obtained from the order $\prec$ by transposing of all $C$-pairs: $c_+\prec_C c$ for all $c\in C$. 
Then, by  for each vertex $v=\{a\prec b\}$ of type \raisebox{-1.0ex}{\m{vgl1}} or \raisebox{-1.0ex}{\m{vgr1}}   there exists a unique adjacent to $v$ vertex $u=\{x\prec y\}$, such that $a\prec_C x \prec_C y \prec_C b$ (we say, that $u$ is nested in $v$). 
We orient this edge from $v$ to $u$ and call it {\it oriented} edge. 
We say that a vertex $v$ of type \raisebox{-1.0ex}{\m{vgl1}} or \raisebox{-1.0ex}{\m{vgr1}} is a {\it head}  if $v$ is not an end point for any oriented edge of $\mathcal{G}$. 
Every head vertex of type \raisebox{-1.0ex}{\m{vgl1}}  continues according to oriented edges by a sequence of vertices  $\raisebox{-1.0ex}{\m{vgr1}}, \raisebox{-1.0ex}{\m{vgl1}},\raisebox{-1.0ex}{\m{vgr1}},...$ ( $\raisebox{-1.0ex}{\m{vgl1}}, \raisebox{-1.0ex}{\m{vgr1}},\raisebox{-1.0ex}{\m{vgl1}},...$ in the case of \raisebox{-1.0ex}{\m{vgr1}}) and ends by a vertex of type \raisebox{-1.0ex}{\m{vgr3}} or \raisebox{-1.0ex}{\m{vgl3}} exactly as it is shown on Fig.~\ref{Fig:RLdef} for $L_k$,$R_l$ for $k,l\ge 3$. 
We call this, subsequent to the head, part of $\mathcal{G}$ a {\it tail} of $\mathcal{G}$.

If a connected component of $\mathcal{G}$ contains a vertex of type \raisebox{-1.0ex}{\m{vgr1}} or \raisebox{-1.0ex}{\m{vgl1}} then it contains at least one head, since there is a  vertex of type \raisebox{-1.0ex}{\m{vgr1}} or \raisebox{-1.0ex}{\m{vgl1}} which is not nested in a sense of $\prec_C$ in any other vertices of type \raisebox{-1.0ex}{\m{vgr1}} or \raisebox{-1.0ex}{\m{vgl1}}. 

\begin{figure}[htb]
\begin{center}
\includegraphics{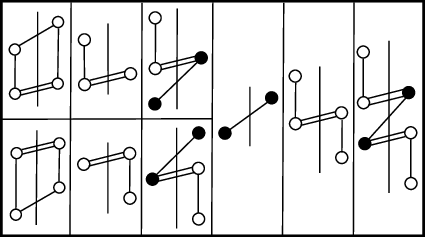}
\caption{\,\, List of all bases without heads}
\label{Fig:Unhead}
\end{center}
\end{figure}

\begin{figure}[htb]
\begin{center}
\includegraphics{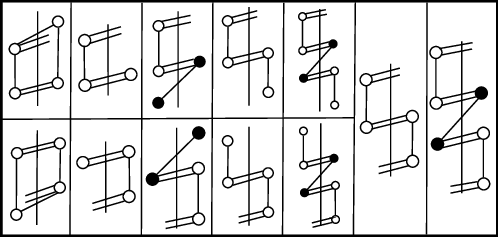}
\caption{\,\, List of all bases without heads}
\label{Fig:Head}
\end{center}
\end{figure}

Consider a connected component of $\mathcal{G}$. It follows from Statement~\ref{val} that it is either homeomorphic to a segment or to a circle, since all valencies of its vertices are 1 or 2.  Let us call a {\it base} of a connected component this component with all tails removed. The list of all bases is finite. Corresponding to bases parts of $\d$ are shown on Fig.~\ref{Fig:Unhead}  and Fig.~\ref{Fig:Head}.

Bases without heads corresponds to direct summands of $\d$ which are equivalent (by action of diagonal subgroup of $\AutT(A,B)$) to summands of type  $L_0\# R_2$, $L_2\# R_0$, $L_1\#\raisebox{-1.5ex}{\m{s2}}$, $\raisebox{-1.5ex}{\m{s3}}\,\,\#R_1$, $L_1\#\raisebox{-1.5ex}{\m{s1}}\#R_0$, $L_0\#\raisebox{-1.5ex}{\m{s1}}\#R_1$, $L_0\#\raisebox{-1.5ex}{\m{s1}}\#R_0$,  $L_1\#R_1$ and $L_1\#\raisebox{-1.5ex}{\m{s1}}\#R_1$ of Theorem~\ref{th1}. Therefore, each component of $\mathcal{G}$ are homotopically trivial. Bases containing a head vertices are shown on  Fig.~\ref{Fig:Head}.

They corresponds  to direct summands of $\d$ which are equivalent (up to an action of diagonal subgroup of $\AutT(A,B)$) to summands $L_0\# R_l (l\ge 3)$, $L_k\# R_0 (k\ge 3)$, $L_k\#\raisebox{-1.5ex}{\m{s2}} (k\ge 2)$, $\raisebox{-1.5ex}{\m{s3}}\,\,\#R_l (l \ge 2)$, $L_k\#\raisebox{-1.5ex}{\m{s1}}\#R_0 (k\ge 2)$, $L_0\#\raisebox{-1.5ex}{\m{s1}}\#R_l (l\ge 2)$, $L_k\#R_1 (k\ge 2)$, $L_1\#R_l (l\ge 2)$, $L_1\#\raisebox{-1.5ex}{\m{s1}}\#R_1$, $L_k\#\raisebox{-1.5ex}{\m{s1}}\#R_1 (k\ge 2)$, $L_1\#\raisebox{-1.5ex}{\m{s1}}\#R_l (l\ge 2)$, $L_k\#R_l (k,l\ge 2)$, $L_k\#\raisebox{-1.5ex}{\m{s1}}\#R_l (k,l~\ge~2)$ of Theorem~\ref{th1}. Now it is easy to check that we get the list of summands of Theorem~\ref{th1}. It completes the proof of Lemma~\ref{summands} and Theorem~\ref{th1}.\qed

\section{An extension of a strong Morse germ} \label{ext}

\subs{Critical points of an extension}

A strong Morse function $F$ on a compact manifold $M$ of dimension $n$ defines a triple of graded linearly ordered sets $A_F\supset B_F\supset C_F$ (see \ref{amodel}).
$C_F$ is a set of outward critical points of $F|_{\lo{\d M}}$.
For each outward critical point $x$ of $F|_{\lo{\d M}}$ the set  $A_F$ contains an element next to $x$.
Also it contains all the critical points of $F$.
A number of critical points of $F$ is equal to $\# A_F-\# B_F-\# C_F$.

The function $F$ defines an $M$-pair $\M_{A_F,B_F,\,\d}$ up to equivalence (see Statement ~\ref{F-CW}).
There is a unique  $M$-pair $M_{A_F,B_F,\,\delta }$ weak equivalent to $\M_{A_F,B_F,\,\d}$ and decomposed  into direct summands of Theorem~\ref{th1}.
For this pair $M_{A_F,B_F,\,\delta }$ all elements of $C_F$ are $\delta$-trivial and a set of its boundary homologically essential elements conside with  boundary homologically essential elements for $\M_{A_F,B_F,\,\d}$.

Critical value $c$ of $F|_{\d M}$ is called boundary homologically essential if a function of a parameter $a$
\[
\dim(\d^*H_*(M,\d M;\E)\cap i_*H_*(\{F|_{\d M}\le a;\E))
\]
jumps at $a=c$.
Here $\d^*$ is a connected homomorphism from long exact sequence of a pair $(M,\d M)$,  a map $i_*$ is induced by a natural inclusion  $\{F|_{\lo{\d M}}\le a\}\hookrightarrow \d M$.
By definition homologically essential critical values of $F|_{\d M}$ depend on $(M, \d M)$ and $F|_{\d M}$ only.
Homologically essential critical points $H(F|_{\d M})$ of $F|_{\d M}$  correspond to homologically essential elements  of $M_{A_F,B_F,\,\delta }$. 

There are a finite number of all linearly ordered graded sets (with values of degree in $\{0,...,n\}$) of cardinality at most $k$.  
There are a finite number of all summands of Theorem~\ref{th1} having total dimension at most $k$.
Hence, there are a finite number of $M_{A,B}$-pairs such that 1) its homologies correspondingly equal to $H_i(M,\d M; \E)$, $H_i(M;\E)$ and $H_i(\d M;\E)$; 2) $B$ is graded ordered isomorphic to the set of all critical points of $g|_{\d M}$ and $\d|_B$ is algebraic model for $g|_{\d M}$ ; 3) all points in $B$ corresponding to outward critical points  are $\d$-trivial;  4) the set of boundary essential elements isomorphic to $H(g|_{\d M})$; 5) $M$-pair is  decomposed into direct sum of summands of Theorem~\ref{th1}.
We say that such an $M$-pair is tamed to $(g,M,\E)$.
Moreover, one can explicitly find smallest $k=\# A$ such that such $M$-pair exists. 
A number $n_1(g,M,\E)=k-\# B_F-\#C_F$ is a (possible) answer to Arnold's question of estimation from below of a number of critical points of  Morse extension of a strong Morse germ.
\begin{theorem}
\label{n_1}
Let $F$ be a Morse function extending $g$ to $M$.
Then total number $T$ of critical points is greater than or equal to $n_1(g,M,\E)$.
\end{theorem}
This result is a generalization of Barannikov result~\cite{Serega} for a ball $M=D^n$.  
It also can be viewed as a generalization of weak Morse estimation  
\[
T\ge \sum \dim H_i(M,\E)
\]
 in the case of closed manifold,  since for the case of empty boundary $n_1(g,M,\E)=\sum \dim H_i(M,\E)$.

Consider now a function $G=-F$. 
It also produces triple of sets $A_G\supset B_G\supset C_G$.
It is clear how this triple of sets connected with $A_F\supset B_F\supset C_F$. 
Each critical point of $F|_{\d M}$ is a critical point for $G|_{\d M}$,  each inward critical point of $F|_{\d M}$ is outward critical point for $G|_{\d M}$ and vice versa, each critical point of $F$ is a critical point for $G$.
It is clear that $A_G$ ($B_G$ correspondingly) is a function of $A_F$ ($B_F$ corr.) and denote it by $\widehat{A_F}$ ($\widehat{B_F}$ corr.).
Consider now a smallest number $l=\# A$ such that there exists a $M$-pair on $(A,B)$  tamed to $(g,M,\E)$ and that there exists a $M$-pair on $(\widehat{A},\widehat{B})$  tamed to $(-g,M,\E)$.
By  $n_2(g,M,\E)$ we denote $l-\# B_F-\#C_F$ .
We get the following result.
\begin{theorem}
\label{n_2}
Let $F$ be a Morse function extending $g$ to $M$.
Then total number of critical points is greater than or equal to $n_2(g,M,\E)$.
\end{theorem}
By construction we get $n_2(g,M,\E)\ge n_1(g,M,\E)$.

\subs{Examples} Consider $\R P^{2n}$.  
We construct a germ of a function along a sphere $S^{2n-1}$ embedded in $\R P^{2n}$.
Let $D^{2n}$ be a (small) coordinate disk in  $\R P^{2n}$,  consider a function  $F=\sum a_ix_i^2+C$ in coordinates on $D^{2n}$,  all $a_i$ are positive ad distinct numbers. 
Consider a sphere $S^{2n - 1}=\{\sum x_i^2=\varepsilon\}=\d D^{2n}(\varepsilon)$.  
A germ of $F$ along $S^{2n - 1}$ is a Morse germ with two critical points of each index $i\in\{0,...,2n-1\}$.
Let $g$ be its $C^\infty$-small perturbation such that $f$ is strong Morse germ.
Consider now a manifold $M=\R P^{2n}\setminus \mathrm{Int} (D^{2n}(\varepsilon))$ with a boundary $S^{2n - 1}$ and a germ $g$ along $S^{2n - 1}$.
Then one can calculate $n_1(g,M,\E)$ and get $n_1(g,M,\E)=2$ for $\mathrm{char}\,\E \ne 2$.  
Surprisingly $n_2(g,M,\E)=2n$,  moreover one can deduce from the construction of $M$-pairs for Theorem~\ref{n_2}  that a number of critical points of index $i$ is at least $1$ for each $i\in\{1,...,2n\}$.
Note than classical Morse estimates~\cite{Morse} cannot guarantee a single point of extension of $g$ to~$M$.

This example is related with the following.
Suppose now that a closed manifold ($\d M=\varnothing$)   $M$ is $\E$-homological point ($H_i(M; \E)=0$ for each $i>0$) and $\dim M=2n$.
Let $\mathrm{char}\,\E \ne 2$.
Such a manifold exists, for example $M=\R P^{2n}$.
Then $n_1(g,M,\E)=1$,  
but $n_2(g,M,\E)=2n+1$. 
Moreover it is easy to show that a number of critical points of index $i$ is at least $1$.

Note that last example is related with Takens and White result \cite{TW},  claiming that a Morse function on a sphere $S^n$ having even number of critical points in each index has greater than or equal to $2n+2$ critical points.

We will consider the Morse inequalities in more detail in subsequent works.

\subs{Criterium for existence of critical point of any extension} 

Estimations of Theorems~\ref{n_1} and \ref{n_2} can guarantee an existence of possible degenerate critical points. of any extension. 
Note that each compact connected manifold with nonempty boundary admits a Morse function without critical points.

\begin{theorem}
\label{criterium}
Suppose $g$ is strong Morse germ along $\d M$ such that $n_2(g,M,\E)\ne 0$ (or $n_1(g,M,\E)\ne 0$)).  
Then any extension of $g$ to the whole $M$ has at least one (possibly degenerate) critical point. 
\end{theorem}
\pr Indeed an extension without critical points inside $M$ is a strong Morse function.  For its algebraic model $k-\# B_F-\#C_F=0.$
\qed

\ign{
\section{Corrections}
\begin{enumerate}
\item Reference to Arnold problem
\item Short intro to Section \ref{formulation}
\item Boundary map or connecting homomorphism?
\item Check the proof of Lemma \ref{dprimitive}.
\item Grading ---> degree
\item Nuzhno li Statement PQRDEF?
\item Ispravit' dokazatel'stwo lemmy 4.* in view of the definition of weakly equivalence
\item Dopisat' v minimal complexes 1. -- mozhno nasil'no opredelit' P po minimal indecomposable. 2. tuda zhe
ili kuda-to esche chto ekwiwalentnost' zadannaya P, est' "usilenie" ili "oslablenie" slaboj i chto klass
P-ekvivalentosti soderzhit konechnoe chislo klassov weak equiv
\item proverit', chto gruppa deistvuet $\d \to S\d S^{-1}$ a ne $\d \to S^{-1}\d S$
\item emptyset -> varnothing
\item grading -> degree proverit' vsyudu
\item $\mathcal{F}(A)$ ili $\E \otimes A$ - (ne) sdelat' (li) edinoobrazno.
\item kak postavit' kvadrat v vyklyuchnoi formule v konec? (Theorem {Morsedecomp} for example)
\item dobavit' ssylku w d-trivial na primer iz topologii i sootw Lemmu w dokazatel'stwe
\item proverit', podumat', ne nado li vsyudu zamenit' $F_a$ na $M^a$ i $F^{\d}_a$ i $\d M^a$
\item underlevel or sublevel
\item Predlozhenie pro algebr model' podvislo
\item znachok izomorphizma $\sim$ ili $approx$??? Proverit' i vsyudu pomenyat'
\item The bibliography - chto nuzhno, chto net.
\end{enumerate}
}


\ign{
\input{mor119intro.tex}

\input{mor219intro.tex}

\input{mor319.tex}

\input{mor419.tex}

\input{mor519.tex}

\input{mor619.tex}

\input{mor1bibl.tex}

}
\end{document}